\newif\ifpreprint
\preprinttrue

\newif\ifmain
\maintrue

\RequirePackage{tikz}
\ifpreprint
    \documentclass{article}
\else
    \documentclass{sn-jnl}
\fi

\usepackage{subfiles}
\usepackage{ifthen}

\newcommand{\TheTitle}{An accelerated proximal gradient method for multiobjective optimization}

\ifpreprint
    \usepackage{authblk}

\title{\TheTitle\footnote{A previous version of this manuscript can be seen in the Department of Applied Mathematics and Physics, Kyoto University's technical report system (\url{http://www.amp.i.kyoto-u.ac.jp/tecrep/ps_file/2022/2022-001.pdf}).}}

\author[1]{Hiroki Tanabe}
\author[2]{Ellen H. Fukuda}
\author[2]{Nobuo Yamashita}

\affil[1]{Yahoo Japan Corporation}
\affil[2]{Kyoto University}
\affil[ ]{\{tanabe.hiroki.45n@kyoto-u.jp\},\{ellen,nobuo\}@i.kyoto-u.ac.jp}

\date{}

\usepackage[totalwidth=480pt, totalheight=680pt]{geometry}

\usepackage{graphicx}

\usepackage{booktabs}
\usepackage{multirow}
\usepackage{tabularray}

\usepackage{algorithm}
\usepackage{algpseudocode}

\usepackage{amsmath,amssymb,amsthm}
\usepackage{cases}
\usepackage{derivative}
\usepackage{mathtools}
\usepackage{hyperref}
\usepackage[capitalize,nameinlink]{cleveref}
\usepackage{autonum}

\usepackage{setspace}
\usepackage{cite}
\usepackage{siunitx}
\usepackage[inline]{enumitem}
\setlist[enumerate,2]{label=(\alph*), ref=(\alph*)}

\def\usebfsetcapital{\def\setcapital##1{\mathbf{##1}}}

\usebfsetcapital

\def\setR{\setcapital{R}}
\def\setRpos{\setcapital{R}_{\mathord{+}}}

\newcommand{\st}{\mathrm{s.t.}}
\newcommand{\level}{\mathcal{L}}

\newcommand\condition[1]{\quad \text{#1}}
\newcommand\forallcondition[1]{\condition{for all~$#1$}}
\newcommand\eqand{\quad \text{and} \quad}
\DeclareMathOperator*{\argmax}{argmax}
\DeclareMathOperator*{\argmin}{argmin}

\DeclareMathOperator{\interior}{int}

\DeclareMathOperator{\dom}{dom}
\DeclareMathOperator{\prox}{\mathbf{prox}}

\DeclareMathOperator{\envelope}{\mathcal{M}}
\DeclareMathOperator{\indicator}{\chi}
\DeclarePairedDelimiter{\abs}{\lvert}{\rvert}
\DeclarePairedDelimiter{\norm}{\lVert}{\rVert}

\DeclarePairedDelimiter{\setof}{\lbrace}{\rbrace}
\DeclarePairedDelimiterX{\Set}[2]{\lbrace}{\rbrace}{#1\mathrel{}\delimsize\vert\mathrel{}#2}
\DeclarePairedDelimiterX{\innerp}[2]{\langle}{\rangle}{#1, #2}

{}

\newcommand{\acc}{\mathrm{acc}}

\newtheorem{theorem}{Theorem}[section]
\newtheorem{lemma}{Lemma}[section]

\newtheorem{proposition}{Proposition}[section]
\newtheorem{assumption}{Assumption}[section]
\newtheorem{corollary}{Corollary}[section]
\newtheorem{remark}{Remark}[section]

\crefname{equation}{}{}
\Crefname{equation}{Eq.}{Eqs.}
\crefname{enumi}{}{}
\crefname{figure}{Figure}{Figures}
\crefname{assumption}{Assumption}{Assumptions}
\crefname{line}{line}{lines}

\setlist[enumerate]{
    label=(\roman*)
}
\AtBeginEnvironment{theorem}{
    \setlist[enumerate]{
        label=(\roman*),
        ref=Theorem~\thetheorem~(\roman*)
    }
}
\AtEndEnvironment{theorem}{
    \setlist[enumerate]{
        label=(\roman*),
        ref=(\roman*)
    }
}
\AtBeginEnvironment{lemma}{
    \setlist[enumerate]{
        label=(\roman*),
        ref=Lemma~\thelemma~(\roman*)
    }
}
\AtEndEnvironment{lemma}{
    \setlist[enumerate]{
        label=(\roman*),
        ref=(\roman*)
    }
}
\AtBeginEnvironment{proposition}{
    \setlist[enumerate]{
        label=(\roman*),
        ref=Proposition~\theproposition~(\roman*)
    }
}
\AtEndEnvironment{proposition}{
    \setlist[enumerate]{
        label=(\roman*),
        ref=(\roman*)
    }
}
\AtBeginEnvironment{corollary}{
    \setlist[enumerate]{
        label=(\roman*),
        ref=Corollary~\thecorollary~(\roman*)
    }
}
\AtEndEnvironment{corollary}{
    \setlist[enumerate]{
        label=(\roman*),
        ref=(\roman*)
    }
}
\AtBeginEnvironment{remark}{
    \setlist[enumerate]{
        label=(\roman*),
        ref=Remark~\theremark~(\roman*)
    }
}
\AtEndEnvironment{remark}{
    \setlist[enumerate]{
        label=(\roman*),
        ref=(\roman*)
    }
}
\AtBeginEnvironment{definition}{
    \setlist[enumerate]{
        label=(\roman*),
        ref=Definition~\thedefinition~(\roman*)
    }
}
\AtEndEnvironment{definition}{
    \setlist[enumerate]{
        label=(\roman*),
        ref=(\roman*)
    }
}
\AtBeginEnvironment{assumption}{
    \setlist[enumerate]{
        label=(\roman*),
        ref=Assumption~\theassumption~(\roman*)
    }
}
\AtEndEnvironment{assumption}{
    \setlist[enumerate]{
        label=(\roman*),
        ref=(\roman*)
    }
}
\AtBeginEnvironment{example}{
    \setlist[enumerate]{
        label=(\roman*),
        ref=Example~\theexample~(\roman*)
    }
}
\AtEndEnvironment{example}{
    \setlist[enumerate]{
        label=(\roman*),
        ref=(\roman*)
    }
}
\newcounter{subcreftmpcnt}
\newcommand\romansubformat[1]{(\roman{#1})}
\makeatletter
\newcommand\subcref[2][\romansubformat]{
\ifcsname r@#2@cref\endcsname
  \cref@getcounter {#2}{\mylabel}
  \setcounter{subcreftmpcnt}{\mylabel}
  \hyperref[#2]{part~#1{subcreftmpcnt}}
 \else ?? \fi}
\newcommand\sublabelcref[2][\romansubformat]{
\ifcsname r@#2@cref\endcsname
  \cref@getcounter {#2}{\mylabel}
  \setcounter{subcreftmpcnt}{\mylabel}
  \hyperref[#2]{#1{subcreftmpcnt}}
 \else ?? \fi}
\newcommand\subCref[2][\romansubformat]{
\ifcsname r@#2@cref\endcsname
  \cref@getcounter {#2}{\mylabel}
  \setcounter{subcreftmpcnt}{\mylabel}
  \hyperref[#2]{Part~#1{subcreftmpcnt}}
 \else ?? \fi}

\makeatother

\else
    \input{shared_coap}
\fi

\begin{document}

\ifpreprint
    \maketitle
\else
\fi

\newcommand\TheAbstract{
    This paper presents an accelerated proximal gradient method for multiobjective optimization, in which each objective function is the sum of a continuously differentiable, convex function and a closed, proper, convex function. Extending first-order methods for multiobjective problems without scalarization has been widely studied, but providing accelerated methods with accurate proofs of convergence rates remains an open problem. Our proposed method is a multiobjective generalization of the accelerated proximal gradient method, also known as the Fast Iterative Shrinkage-Thresholding Algorithm (FISTA), for scalar optimization. The key to this successful extension is solving a subproblem with terms exclusive to the multiobjective case. This approach allows us to demonstrate the global convergence rate of the proposed method~($O(1 / k^2)$), using a merit function to measure the complexity. Furthermore, we present an efficient way to solve the subproblem via its dual representation, and we confirm the validity of the proposed method through some numerical experiments.
}
\ifpreprint
    \begin{abstract}
        \TheAbstract
    \end{abstract}
\fi

\section{Introduction} \label{sec: intro}
Multiobjective optimization consists in minimizing (or maximizing) more than one objective function at once under possible constraints.
In general, there is no single point that minimizes all objective functions simultaneously, so the concept of \emph{Pareto optimality} becomes essential.
We call a point Pareto optimal if there is no other point with the same or smaller objective function values and with at least one objective function value being strictly smaller.

One of the most popular strategies for solving multiobjective optimization problems is the \emph{scalarization approach}~\cite{Gass1955,Geoffrion1968,Zadeh1963}.
It converts the original multiobjective problem into another, which has a parametrized scalar-valued objective function.
If each objective function in the multiobjective optimization problem is convex, the converted single objective optimization problems typically become convex optimization.
However, it can be challenging to choose the appropriate parameters (or weights) in advance.
For example, Marler and Arora~\cite{Marler2010} discussed such difficulties in the weighted sum method.
Another approach, which does not use scalarization, is based on metaheuristics~\cite{Gandibleux2004} but lacks a theoretical proof of convergence to Pareto solutions.

To overcome those drawbacks, many descent algorithms for multiobjective optimization problems have been developed recently~\cite{Fukuda2014}.
These algorithms decrease all objective functions at each iteration, offer the advantages of not requiring a priori parameter selection, and provide convergence guarantees under reasonable assumptions.
For instance, Fliege and Svaiter~\cite{Fliege2000} proposed the steepest descent method for differentiable multiobjective optimization problems.
Other examples include the projected gradient~\cite{Fukuda2013,Grana-Drummond2004,Zhao2021}, Newton's~\cite{Fliege2009,Goncalves2022}, trust-region~\cite{Carrizo2016}, and conjugate gradient methods~\cite{LucambioPerez2018}.
Descent methods for infinite-dimensional vector optimization problems have also been studied, including the proximal point~\cite{Bonnel2005} and the inertial forward-backward methods~\cite{Bot2018}.
Among these, methods that use only the first-order derivatives of the objective functions, such as the steepest descent and the projected gradient methods, are called \emph{first-order methods}.
Another well-known multiobjective first-order method is the proximal gradient~\cite{Tanabe2019}, which works for composite problems, i.e., with each objective being the sum of a differentiable function and a convex but not necessarily differentiable one.
This algorithm, as well as the steepest descent, is known to converge to Pareto solutions with rate~$O(1/k)$~\cite{Fliege2019, Tanabe2023}.

On the other hand, there are many studies related to the acceleration of single-objective first-order methods.
After being established by Nesterov~\cite{Nesterov1983}, researchers developed various accelerated schemes.
In particular, the Fast Iterative Shrinkage-Thresholding Algorithm (FISTA)~\cite{Beck2009}, an accelerated version of the proximal gradient method, has contributed to a wide range of research fields, including image and signal processing.
These methods may increase the objective function values in some iterations, but overall they are known to converge faster than the original descent methods, both theoretically and experimentally.

However, in the multiobjective case, the studies associated with accelerated algorithms are still insufficient~\cite{ElMoudden2021,Wang2020}.
In 2020, El Moudden and El Mouatasim~\cite{ElMoudden2021} proposed an accelerated diagonal steepest descent method for multiobjective optimization, a natural extension of Nesterov's accelerated method for single-objective problems.
They proved the global convergence rate of the algorithm ($O(1 / k^2)$) under the assumption that the sequence of the Lagrange multipliers of the subproblems is eventually fixed.
Nevertheless, this assumption is restrictive because it indicates that the approach is essentially the same as the (single-objective) Nesterov's method, only applied to the minimization of a weighted sum of the objective functions.

Here, we propose a genuine accelerated proximal gradient method for multiobjective optimization.
As it is usual, in each iteration, we solve a convex (scalar-valued) subproblem.
While the accelerated and non-accelerated algorithms solve the same subproblem in the single-objective case, the subproblem of our accelerated method has terms that are not included in the non-accelerated version.
However, we can ignore these terms in the single-objective case, and thus we can regard our proposed method as a generalization of FISTA.
Moreover, under more natural assumptions, we prove the proposed method's global convergence rate ($O(1/k^2)$) by using a merit function~\cite{Tanabe2023a} to measure the complexity.

Furthermore, having the practical computational efficiency in mind, we derive a dual of the subproblem, which is convex and differentiable.
Such a dual problem turns out to be easier to solve than the original one, especially when the number of objective functions is smaller than the dimension of the decision variables.
We can also reconstruct the original subproblem's solution directly from the dual optimum.
In addition, we implement the whole algorithm using this dual problem and confirm its effectiveness with numerical experiments.

The outline of this paper is as follows.
In \cref{sec: preliminaries}, we introduce some notations and concepts used in this paper.
\Cref{sec: prox} recalls the proximal gradient method for multiobjective optimization proposed in~\cite{Tanabe2019}.
We present the proposed accelerated proximal gradient method for multiobjective optimization in \cref{sec: acc prox} and analyze its~$O(1 / k^2)$ convergence rate in \cref{sec: convergence rate}.
Moreover, \cref{sec: subproblem} introduces an efficient way to solve the subproblem via its dual form.
Finally, we report some numerical results for test problems in \cref{sec: experiments}, demonstrating that the proposed method is faster than the one without acceleration.

\ifmain
\else
\fi

\section{Preliminaries} \label{sec: preliminaries}
All over this work, for any natural number~$d$,~$\setR^d$ denotes the $d$-dimensional real space, $\setRpos^d \subseteq \setR^d$ designates the nonnegative orthant of~$\setR^d$, i.e.,$\setRpos^d \coloneqq \Set*{v \in \setR^d}{v_i \ge 0, i = 1, \dots, d}$, and~$\Delta^d$ represents the standard simplex in~$\setR^d$ given by
\[ \label{eq: simplex}
    \Delta^d \coloneqq \Set*{\lambda \in \setRpos^d}{\sum_{i = 1}^{d} \lambda_i = 1}
.\] 
Then, we can consider the partial orders induced by~$\setRpos^d$: for all~$v^1, v^2 \in \setR^d$,~$v^1 \le v^2$ (alternatively, $v^2 \ge v^1$) if~$v^2 - v^1 \in \setRpos^d$ and~$v^1 < v^2$ (alternatively,~$v^2 > v^1$) if~$v^2 - v^1 \in \interior \setRpos^d$.
In other words, $v^1 \le v^2$ and~$v^1 < v^2$ stand for~$v^1_i \le v^2_i$ and~$v^1_i < v^2_i$ for all~$i = 1, \dots, d$, respectively.
Moreover, let~$\innerp*{\cdot}{\cdot}$ be the Euclidean inner product in~$\setR^d$, i.e.,~$\innerp*{u}{v} \coloneqq \sum_{i = 1}^d {u_i v_i}$, and let~$\norm*{ \cdot }$ be the Euclidean norm, i.e., $\norm*{u} \coloneqq \sqrt{\innerp*{u}{u}}$.
Furthermore, we define the~$\ell_1$-norm and the~$\ell_\infty$-norm by~$\norm*{u}_1 \coloneqq \sum_{i = 1}^{d} \abs*{u_i}$ and~$\norm*{u}_\infty \coloneqq \max_{i = 1, \dots, d} \abs*{u_i}$, respectively.
We now recall the obvious equality related to norm and inner product:
\[ \label{eq: Pythagoras}
    \norm{b - a}^2 + 2 \innerp{b - a}{a - c} = \norm{ b - c }^2 - \norm{a - c}^2 \forallcondition{a, b, c \in \setR^n}
.\] 

On the other hand, for a closed, proper and convex function~$h \colon \setR^n \to \setR \cup \setof*{+ \infty}$, we call~$\eta \in \setR^n$ a subgradient of~$h$ at~$x \in \setR^n$ if
\[ \label{eq: subgradient}
    h(y) \ge h(x) + \innerp*{\eta}{y - x} \forallcondition{y \in \setR^n}
,\]
and we write~$\partial h(x)$ the subdifferential of~$h$ at~$x$, i.e., the set of all subgradients of~$h$ at~$x$.
In addition, the subdifferential for a vector-valued function is the direct product of the subdifferentials of each component.
We also define the \emph{Moreau envelope} or \emph{Moreau-Yosida regularization}~\cite{Moreau1965,Yosida1995} of~$h$ by
\[ \label{eq: Moreau envelope}
    \envelope_h(x) \coloneqq \min_{y \in \setR^n} \left[ h(y) + \frac{1}{2} \norm*{x - y}^2 \right] 
.\] 
The minimization problem in~\cref{eq: Moreau envelope} has a unique solution because of the strong convexity of its objective function.
We call this solution the \emph{proximal operator} and write it as
\[ \label{eq: proximal operator}
    \prox_h(x) \coloneqq \argmin_{y \in \setR^n} \left[ h(y) + \frac{1}{2} \norm*{x - y}^2 \right] 
.\] 
\begin{remark} \label{rem: Moreau and prox}
    \begin{enumerate}
        \item \cite[Theorem 6.24]{Beck2017} If~$h$ is the \emph{indicator function} of a nonempty set~$S \subseteq \setR^n$, i.e.,
            \[ \label{eq: indicator}
                \indicator_S(x) = \begin{dcases}
                    0, & x \in S, \\
                    + \infty, & x \notin S
                ,\end{dcases}
            \] 
            then the proximal operator reduces to the projection onto~$S$. \label{enum: prox of indicator}
        \item \cite[Theorem 6.42]{Beck2017} The proximal operator of a closed, proper, and convex function~$h$ is non-expansive, i.e.,~$\norm*{\prox_h(x) - \prox_h(y)} \le \norm*{x - y}$.
            In other words, $\prox_h$ is~$1$-Lipschitz continuous.
        \item \cite[Theorem 6.60]{Beck2017} Even if a closed, proper, and convex function~$h$ is non-differentiable, its Moreau envelope~$\envelope_h$ has a~$1$-Lipschitz continuous gradient as follows:~$\nabla \envelope_h(x) = x - \prox_h(x)$.
    \end{enumerate}
\end{remark}

We now focus on the following multiobjective optimization problem:
\[ \label{eq: MOP}
    \min_{x \in \setR^n} \quad F(x)
\]
with a vector-valued function~$F \colon \setR^n \to (\setR \cup \setof{+ \infty})^m$ with~$F \coloneqq (F_1, \dots, F_m)^\top$.
We assume that each component~$F_i \colon \setR^n \to \setR \cup \setof{+ \infty}$ is defined by~$F_i(x) \coloneqq f_i(x) + g_i(x)$ for all~$i = 1, \dots, m$ with convex and continuously differentiable functions~$f_i \colon \setR^n \to \setR, i = 1, \dots, m$ and closed, proper and convex functions~$g_i \colon \setR^n \to \setR \cup \setof*{+ \infty}, i = 1, \dots, m$.
We also suppose that each~$\nabla f_i$ is Lipschitz continuous with constant~$L_i > 0$ and define~$L \coloneqq \max_{i = 1, \dots, m} {L_i}$.
From the so-called descent lemma~\cite[Proposition A.24]{Bertsekas1999}, we have~$f_i(p) - f_i(q) \le \innerp*{\nabla f_i(q)}{p - q} + (L / 2) \norm*{p - q}^2$
for all~$p, q \in \setR^n$ and~$i = 1, \dots, m$, which gives
\[ \label{eq: descent mid}
    \begin{aligned}
        F_i(p) - F_i(r) &= f_i(p) - f_i(q) + g_i(p) + f_i(q) - F_i(r) \\
                        &\le \innerp*{\nabla f_i(q)}{p - q} + g_i(p) + f_i(q) - F_i(r) + \frac{L}{2} \norm*{p - q}^2
    \end{aligned}
\]
for all~$p, q, r \in \setR^n$ and~$i = 1, \dots, m$.

Now, we introduce some concepts used in the multiobjective optimization problem~\cref{eq: MOP}.
Recall that
\[ \label{eq: weak Pareto}
    X^\ast \coloneqq \Set*{x^\ast \in \setR^n}{\text{There does not exist~$x \in \setR^n$ such that~$F(x) < F(x^\ast)$}}
\]
is the set of \emph{weakly Pareto optimal} points for~\cref{eq: MOP}.
We also define the effective domain of~$F$ by~$\dom F \coloneqq \Set*{x \in \setR^n}{F_i(x) < + \infty \text{ for all } i = 1, \dots, m}$,
and we write the level set of~$F$ on~$c \in \setR^m$ as
\[ \label{eq: level set}
    \level_F(c) \coloneqq \Set*{x \in \setR^n}{F(x) \le c}
.\]
In addition, we express the image of~$A \subseteq \setR^n$ and the inverse image of~$B \subseteq (\setR \cup \setof*{+ \infty})^m$ under~$F$ as~$F(A) \coloneqq \Set*{F(x) \in \setR^m}{x \in A} \eqand F^{-1}(B) \coloneqq \Set*{x \in \setR^n}{F(x) \in B}$,
respectively.

Finally, let us recall the merit function~$u_0 \colon \setR^n \to \setR \cup \setof*{+ \infty}$ proposed in~\cite{Tanabe2023a}:
\[ \label{eq: u_0}
    u_0(x) \coloneqq \sup_{z \in \setR^n} \min_{i = 1, \dots, m} [ F_i(x) - F_i(z) ]
,\]
which returns zero at optimal solutions and strictly positive values otherwise.
The following theorem shows that $u_0$ is a merit function in the Pareto sense.
\begin{theorem}~\cite[Theorem 3.1]{Tanabe2023a} \label{thm: merit Pareto}
    Let~$u_0$ be defined by~\cref{eq: u_0}.
    Then, we get~$u_0(x) \ge 0$ for all~$x \in \setR^n$.
    Moreover,~$x \in \setR^n$ is weakly Pareto optimal for~\cref{eq: MOP} if and only if~$u_0(x) = 0$.
\end{theorem}
Note that when~$m = 1$, we have~$u_0(x) = F_1(x) - F_1^\ast$, where~$F_1^\ast$ is the optimal objective value.
This is clearly a merit function for scalar-valued optimization.

\ifmain
\else
\fi

\section{Proximal gradient methods for multiobjective optimization} \label{sec: prox}
Let us now recall the proximal gradient method for~\cref{eq: MOP}, an extension of the classical proximal gradient method, proposed by Tanabe, Fukuda, and Yamashita~\cite{Tanabe2019}.
We explain how to generate the sequence of iterates, and afterward, we show the algorithm and its convergence rate.

For given~$x \in \dom F$ and~$\ell > 0$, we consider the following minimization problem:
\[ \label{eq: prox subprob}
    \min_{z \in \setR^n} \quad \varphi_\ell(z; x)
,\]
where~$\varphi_\ell(z; x) \coloneqq \max_{i = 1, \dots, m} [ \innerp*{\nabla f_i(x)}{z - x} + g_i(z) - g_i(x) ] + (\ell / 2) \norm{ z - x }^2$.
The convexity of~$g_i$ implies that~$z \mapsto \varphi_\ell(z; x)$ is strongly convex, so the problem~\cref{eq: prox subprob} always has a unique solution.
Let us write such a solution as~$p_\ell(x)$ and let~$\theta_\ell(x)$ be its optimal function value, i.e.,
\[ \label{eq: p theta}
    p_\ell(x) \coloneqq \argmin_{z \in \setR^n} {\varphi_\ell(z; x)} \eqand \theta_\ell(x) \coloneqq \min_{z \in \setR^n} \varphi_\ell(z; x)
.\]
The following proposition shows that~$p_\ell(x)$ and~$\theta_\ell(x)$ helps to characterize the weak Pareto optimality of~\cref{eq: MOP}.
\begin{proposition} \label{prop: prox terminate}
    Let~$p_\ell$ and~$\theta_\ell$ be defined by~\cref{eq: p theta}.
    Then, the statements below hold.
    \begin{enumerate}
        \item The following three conditions are equivalent: 
            \begin{enumerate*}
                \item $x$ is weakly Pareto optimal; 
                \item $p_\ell(x) = x$;
                \item $\theta_\ell(x) = 0$.
            \end{enumerate*}
        \item The mappings~$p_\ell$ and~$\theta_\ell$ are both continuous.
    \end{enumerate}
\end{proposition}
\begin{proof}
    It is clear from~\cite[Lemma 3.2]{Tanabe2019} and the convexity of~$f_i$.
\end{proof}

From \cref{prop: prox terminate}, we can treat~$\norm*{p_\ell(x) - x}_\infty < \varepsilon$ for some~$\varepsilon > 0$ as a stopping criteria.
Moreover, if~$\ell > L / 2$ then we have~$F_i(p_\ell(x)) \le F_i(x)$ for all~$x \in \dom F$ and~$i = 1, \dots, m$~\cite{Tanabe2023}.
Now, we state below the proximal gradient method for~\cref{eq: MOP}.

\begin{algorithm}[hbtp]
    \caption{Proximal gradient method for multiobjective optimization~\cite{Tanabe2019}}
    \label{alg: pgm}
    \begin{algorithmic}[1]
        \Require $x^0 \in \dom F$, $\ell > L / 2$, $\varepsilon > 0$
        \Ensure $x^\ast$: A weakly Pareto optimal point
        \State $k \gets 0$
        \Loop
        \State $x^{k + 1} \gets p_\ell(x^k)$, where~$p_\ell$ is defined by~\cref{eq: p theta}
        \If{$\norm*{x^{k + 1} - x^k}_\infty \ge \varepsilon$}
        \Return $x^{k + 1}$
        \EndIf
        \State $k \gets k + 1$
        \EndLoop
    \end{algorithmic}
\end{algorithm}

When~$\ell \ge L$, \cref{alg: pgm} is known to generate~$\setof{x^k}$ such that~$\setof{u_0(x^k)}$ converges to zero with rate~$O(1 / k)$ under the following assumption.
Note that this assumption is not particularly strong, as suggested in~\cite[Remark 5.2]{Tanabe2023}.
\begin{assumption}~\cite[Assumption 5.1]{Tanabe2023} \label{asm: bound}
    Let~$X^\ast$ and~$\level_F$ be defined by~\cref{eq: weak Pareto,eq: level set}, respectively.
    Then, for all~$x \in \level_F(F(x^0))$, there exists~$x^\ast \in X^\ast$ such that~$F(x^\ast) \le F(x)$ and
    \[ \label{eq: R}
        R \coloneqq \sup_{F^\ast \in F(X^\ast \cap \level_F(F(x^0)))} \inf_{z \in F^{-1}(\setof*{F^\ast})} \norm*{z - x^0}^2 < + \infty
    .\]
\end{assumption}
\begin{theorem}~\cite[Theorem 5.2]{Tanabe2023} \label{thm: conv rate}
    Assume that~$\ell \ge L$.
    Then, under \cref{asm: bound}, \cref{alg: pgm} generates a sequence~$\setof{x^k}$ such that~$u_0(x^k) \le (\ell R) / (2 k)$ for all~$k \ge 1$.
\end{theorem}

At the end of this section, we note some remarks about \cref{alg: pgm}.
\begin{remark} \label{rem: proximal gradient method}
    \begin{enumerate}
        \item Since~$x \in \dom F$ implies~$p_\ell(x) \in \dom F$, \cref{alg: pgm} is well-defined.
        \item If~$g_i = 0$, \cref{alg: pgm} corresponds to the steepest descent method~\cite{Fliege2000}:
            \[
                x^{k + 1} \coloneqq \min_{z \in \setR^n} \left[ \max_{i = 1, \dots, m} \innerp*{\nabla f_i(x^k)}{z - x^k} + \frac{\ell}{2} \norm*{z - x^k}^2 \right] 
            .\] 
            On the other hand, when~$f_i = 0$, it matches the proximal point method~\cite{Bonnel2005}:
            \[
                x^{k + 1} \coloneqq \min_{z \in \setR^n} \left\{ \max_{i = 1, \dots, m} \left[ g_i(z) - g_i(x) \right]  + \frac{\ell}{2} \norm*{z - x^k}^2 \right\}
            .\] 
            Furthermore, when~$g_i$ is the indicator function~\cref{eq: indicator} of a convex set~$S \subseteq \setR^n$, it coincides with the projected gradient method~\cite{Grana-Drummond2004}:
            \[
                x^{k + 1} \coloneqq \min_{z \in S - x^k} \left[ \max_{i = 1, \dots, m} \innerp*{\nabla f_i(x^k)}{z - x^k} + \frac{\ell}{2} \norm*{z - x^k}^2 \right] 
            .\] \label{enum: other methods}
        \item When it is difficult to estimate the Lipschitz constant~$L$, we can set the initial value of~$\ell$ appropriately.
            Then, at each iteration we increase~$\ell$ by multiplying it with some prespecified scalar, until~$F_i(p_\ell(x^k)) - F_i(x^k) \le \theta_\ell(x^k)$ is satisfied for all~$i = 1, \dots, m$.
    If~$L$ is finite, the number of times that~$\ell$ is increased is at most a constant.
    \label{enum: backtracking}
    \end{enumerate}
\end{remark}

\ifmain
\else
\fi

\section{An accelerated proximal gradient method for multiobjective optimization} \label{sec: acc prox}
This section proposes an accelerated version of the proximal gradient method for multiobjective optimization.
Similarly to the non-accelerated version given in the last section, a subproblem is considered in each iteration.
More specifically, the proposed method solves the following subproblem for given~$x \in \dom F$,~$y \in \setR^n$, and~$\ell \ge L$:
\[ \label{eq: acc prox subprob}
    \min_{z \in \setR^n} \quad \varphi^\acc_\ell(z; x, y) 
,\]
where
\[ \label{eq: varphi acc}
    \varphi^\acc_\ell(z; x, y) \coloneqq \max_{i = 1, \dots, m} \left[ \innerp{\nabla f_i(y)}{z - y} + g_i(z) + f_i(y) - F_i(x) \right] + \frac{\ell}{2} \norm{z - y}^2
.\]
Note that when~$y = x$,~\cref{eq: acc prox subprob} is reduced to the subproblem~\cref{eq: prox subprob} of the proximal gradient method.
Note also that when~$m = 1$, the subproblem becomes
\[ \label{eq: single}
    \min_{z \in \setR^n} \quad \innerp*{\nabla f_1(y)}{z - y} + g_1(z) + \frac{\ell}{2} \norm{z - y}^2
,\]
which is the subproblem of the single-objective FISTA~\cite{Beck2009}.
The distinctive feature of our proposal~\cref{eq: acc prox subprob} is the term~$f_i(y) - F_i(x)$, whereas the easy analogy from the single-objective subproblem~\cref{eq: single} is
\[ \label{eq: analogy}
    \min_{z \in \setR^n} \quad \max_{i = 1, \dots, m} \left[ \innerp{\nabla f_i(y)}{z - y} + g_i(z) \right] + \frac{\ell}{2} \norm{z - y}^2
.\] 
By putting such a term, the inside of the~$\max$ operator approximates~$F_i(z) - F_i(x)$ rather than~$F_i(z) - F_i(y)$.
This is a negligible difference in the single-objective case, but deeply affects the proof in the multi-objective case.

Since~$g_i$ is convex for all~$i = 1, \dots, m$, $z \mapsto \varphi^\acc_\ell(z; x, y)$ is strongly convex.
Thus, the subproblem~\cref{eq: acc prox subprob} has a unique optimal solution~$p^\acc_\ell(x, y)$ and takes the optimal function value~$\theta^\acc_\ell(x, y)$, i.e.,
\[ \label{eq: p theta acc}
    p^\acc_\ell(x, y) \coloneqq \argmin_{z \in \setR^n} \varphi^\acc_\ell(z; x, y) \eqand \theta^\acc_\ell(x, y) \coloneqq \min_{z \in \setR^n} \varphi^\acc_\ell(z; x, y)
.\]
Moreover, the optimality condition of~\cref{eq: acc prox subprob} implies that for all $x \in \dom F$ and~$y \in \setR^n$ there exists~$\eta(x, y) \in \partial g(p^\acc_\ell(x, y))$ and a Lagrange multiplier~$\lambda(x, y) \in \setR^m$ such that
\begin{subequations} \label{eq: kkt}
    \begin{gather} 
        \sum_{i = 1}^m \lambda_i(x, y) \left[ \nabla f_i(y) + \eta_i(x, y) \right] = - \ell \left[ p^\acc_\ell(x, y) - y \right] \label{eq: optimal} \\
        \lambda(x, y) \in \Delta^m, \quad \lambda_j(x, y) = 0 \forallcondition{j \notin \mathcal{I}(x, y)} \label{eq: lambda}
    ,\end{gather}
\end{subequations}
where~$\Delta^m$ denotes the standard simplex~\cref{eq: simplex} and
\[ \label{eq: I}
    \mathcal{I}(x, y) \coloneqq \argmax_{i = 1, \dots, m} \left[ \innerp{\nabla f_i(y)}{p^\acc_\ell(x, y) - y} + g_i(p^\acc_\ell(x, y)) + f_i(y) - F_i(x) \right]
.\]
Now, we introduce a relation useful for the subsequent analysis.
\begin{lemma} \label{thm: useful relation}
    Let~$p^\acc_\ell$ and~$\theta^\acc_\ell$ be defined by~\cref{eq: p theta acc}.
    Then, we have
    \begin{align}
        \MoveEqLeft - \frac{\ell}{2} \left[ \norm{p^\acc_\ell(x, y) - z}^2 - \norm{y - z}^2 \right] \\
        &\ge \theta^\acc_\ell(x, y) + \sum_{i = 1}^{m} \lambda_i(x, y) \left[ \innerp{\nabla f_i(y)}{y - z} - g_i(z) - f_i(y) + F_i(x) \right] 
    \end{align}
    for all~$x, z \in \dom F$ and~$y \in \setR^n$.
\end{lemma}
\begin{proof}
    Let~$x, z \in \dom F$ and~$y \in \setR^n$.
    From~\cref{eq: optimal} and the definition~\cref{eq: subgradient} of the subgradient, we get
    \begin{align}
        \MoveEqLeft - \ell \innerp{p^\acc_\ell(x, y) - y}{p^\acc_\ell(x, y) - z} \\
        \ge{}& \sum_{i = 1}^{m} \lambda_i(x, y) \left[ \innerp{\nabla f_i(y)}{p^\acc_\ell(x, y) - z} + g_i(p^\acc_\ell(x, y)) - g_i(z) \right] \\
        ={}& \sum_{i = 1}^{m} \lambda_i(x, y) \left[ \innerp{\nabla f_i(y)}{p^\acc_\ell(x, y) - y} + g_i(p^\acc_\ell(x, y)) + f_i(y) - F_i(x) \right] \\
         &+ \sum_{i = 1}^{m} \lambda_i(x, y) \left[ \innerp{\nabla f_i(y)}{y - z} - g_i(z) - f_i(y) + F_i(x) \right] \\
        ={}& \max_{i = 1, \dots, m} \left[ \innerp{\nabla f_i(y)}{p^\acc_\ell(x, y) - y} + g_i(p^\acc_\ell(x, y)) + f_i(y) - F_i(x) \right] \\
         &+ \sum_{i = 1}^{m} \lambda_i(x, y) \left[ \innerp{\nabla f_i(y)}{y - z} - g_i(z) - f_i(y) + F_i(x) \right]
    ,\end{align}
    where the second equality comes from~\cref{eq: lambda,eq: I}.
    Adding~$(\ell / 2) \norm{p^\acc_\ell(x, y) - y}^2$ to both sides and the definition~\cref{eq: p theta acc} of~$p^\acc_\ell$ and~$\theta^\acc_\ell$ lead to
    \begin{align}
        \MoveEqLeft - \frac{\ell}{2} \left[ 2 \innerp{p^\acc_\ell(x, y) - y}{p^\acc_\ell(x, y) - z} - \norm{p^\acc_\ell(x, y) - y}^2 \right] \\
        &\ge \theta^\acc_\ell(x, y) + \sum_{i = 1}^{m} \lambda_i(x, y) \left[ \innerp{\nabla f_i(y)}{y - z} - g_i(z) - f_i(y) + F_i(x) \right]
    .\end{align}
    The left-hand side of this inequality is equal to~$- (\ell / 2) [ 2 \innerp{p^\acc_\ell(x, y) - y}{y - z} + \norm{p^\acc_\ell(x, y) - y}^2 ]$.
    Hence, applying~\cref{eq: Pythagoras} with~$(a, b, c) \coloneqq (y, z, p^\acc_\ell(x, y))$, we get the desired inequality.
\end{proof}
We also note that by taking~$z = y$ in the objective function of~\cref{eq: acc prox subprob}, we have
\[ \label{eq: acc prox subprob optimality}
    \theta^\acc_\ell(x, y) \le \varphi^\acc_\ell(y; x, y) = \max_{i = 1, \dots, m}\left\{ F_i(y) - F_i(x) \right\}
\]
for all~$x \in \dom F$ and~$y \in \setR^n$.
Moreover, from~\cref{eq: descent mid} with~$p = z, q = y, r = x$, and the fact that~$\ell \ge L$, it follows that
\[ \label{eq: acc approx ineq}
    \theta^\acc_\ell(x, y) \ge \max_{i = 1, \dots, m}\left\{ F_i(p^\acc_\ell(x, y)) - F_i(x) \right\}
\]
for all~$x \in \dom F$ and~$y \in \setR^n$.
We now characterize weak Pareto optimality in terms of the mappings~$p^\acc_\ell$ and $\theta^\acc_\ell$, similarly to \cref{prop: prox terminate} for the proximal gradient method.
\begin{proposition} \label{prop: acc prox termination}
    Let~$p^\acc_\ell(x, y)$ and~$\theta^\acc_\ell(x, y)$ be defined by~\cref{eq: p theta acc}.
    Then, the statements below hold.
    \begin{enumerate}
        \item The following three conditions are equivalent: 
            \begin{enumerate*}
                \item $y \in \setR^n$ is weakly Pareto optimal for~\cref{eq: MOP};
                \item $p^\acc_\ell(x, y) = y$ for some~$x \in \setR^n$;
                \item $\theta^\acc_\ell(x, y) = \max_{i = 1, \dots, m} [ F_i(y) - F_i(x) ]$ for some~$x \in \setR^n$.
            \end{enumerate*}\label{enum: acc prox optimality}
        \item The mappings~$p^\acc_\ell$ and~$\theta^\acc_\ell$ are locally H\"{o}lder continuous with exponent~$1 / 2$ and locally Lipschitz continuous, respectively, i.e., for any bounded set~$W \subseteq \setR^n$, there exists~$M_p > 0$ and~$M_\theta > 0$ such that
            \begin{align}
                \norm{p^\acc_\ell(\hat{x}, \hat{y}) - p^\acc_\ell(\check{x}, \check{y})} &\le M_p \norm{(\hat{x}, \hat{y}) - (\check{x}, \check{y})}^{1 / 2}, \\
                \abs{\theta^\acc_\ell(\hat{x}, \hat{y}) - \theta^\acc_\ell(\check{x}, \check{y})} &\le  M_\theta \norm{(\hat{x}, \hat{y}) - (\check{x}, \check{y})}
            \end{align}
            for all~$\hat{x}, \hat{y}, \check{x}, \check{y} \in W$. \label{enum: continuity}
    \end{enumerate}
\end{proposition}
\begin{proof}
    \subCref{enum: acc prox optimality}: From~\cref{eq: acc prox subprob optimality} and the fact that~$\theta^\acc_\ell(x, y) = \varphi^\acc_\ell(p^\acc_\ell(x, y); x, y)$, the equivalence between~(b) and~(c) is apparent.
    Now, let us show that~(a) and~(b) are equivalent.
    When~$y$ is weakly Pareto optimal, we can immediately see from \cref{prop: prox terminate} that~$p^\acc_\ell(x, y) = p_\ell(y) = y$ by letting~$x = y$.
    Conversely, suppose that~$p^\acc_\ell(x, y) = y$ for some~$x \in \setR^n$.
    Let~$z \in \setR^n$ and~$\alpha \in (0, 1)$.
    The optimality of~$p^\acc_\ell(x, y) = y$ for~\cref{eq: acc prox subprob} gives
    \begin{multline}
        \max_{i = 1, \dots, m} \left[ F_i(y) - F_i(x) \right] \le \varphi^\acc_\ell(y + \alpha (z - y); x, y) \\
        = \max_{i = 1, \dots, m} \left[ \innerp*{\nabla f_i(y)}{\alpha (z - y)} + g_i(y + \alpha (z - y)) + f_i(y) - F_i(x) \right] \\
        + \frac{\ell}{2} \norm*{\alpha (z - y)}^2
    .\end{multline}
    Thus, from the convexity of~$f_i$, we get
    \[
            \max_{i = 1, \dots, m}\left[ F_i(y) - F_i(x) \right] \le \max_{i = 1, \dots, m}\left[ F_i(y + \alpha (z - y)) - F_i(x) \right] + \frac{\ell}{2} \norm*{\alpha (z - y)}^2
    .\]
    Moreover, the convexity of~$F_i$ yields
    \begin{multline}
        \max_{i = 1, \dots, m}\left[ F_i(y) - F_i(x) \right] \\
        \begin{aligned}
        &\le \max_{i = 1, \dots, m}\left[  \alpha F_i(z) + (1 - \alpha) F_i(y) - F_i(x)  \right] + \frac{\ell}{2} \norm*{\alpha (z - y) }^2 \\
        &\le \alpha \max_{i = 1, \dots, m}\left[  F_i(z) - F_i(y) \right] + \max_{i = 1, \dots, m}\left\{ F_i(y) - F_i(x)  \right\} + \frac{\ell}{2} \norm*{\alpha (z - y) }^2
        .\end{aligned}
    \end{multline}
    Therefore, we get
    \[
        \max_{i = 1, \dots, m}[  F_i(z) - F_i(y)  ] \ge - \frac{\ell \alpha}{2} \norm*{z - y }^2
    .\]
    Taking~$\alpha \searrow 0$, we obtain~$\max_{i = 1, \dots, m}[  F_i(z) - F_i(y)  ] \ge 0$,
    which implies the weak Pareto optimality of~$y$.

    \subCref{enum: continuity}:
    Take~$\hat{x}, \hat{y}, \check{x}, \check{y} \in W$.
    Adding the two inequalities of~\cref{thm: useful relation} with~$(x, y, z) \coloneqq (\hat{x}, \hat{y}, p^\acc_\ell(\check{x}, \check{y})), (\check{x}, \check{y}, p^\acc_\ell(\hat{x}, \hat{y}))$ gives
    \begin{align}
        \MoveEqLeft - \ell \norm{p^\acc_\ell(\hat{x}, \hat{y}) - p^\acc_\ell(\check{x}, \check{y})}^2 + \frac{\ell}{2} \norm{p^\acc_\ell(\check{x}, \check{y}) - \hat{y}}^2 + \frac{\ell}{2} \norm{p^\acc_\ell(\hat{x}, \hat{y}) - \check{y}}^2 \\
        \ge{}& \theta^\acc_\ell(\hat{x}, \hat{y}) + \theta^\acc_\ell(\check{x}, \check{y}) \\
           &+ \sum_{i = 1}^{m} \lambda_i(\hat{x}, \hat{y}) \left[ \innerp{\nabla f_i(\hat{y})}{\hat{y} - p^\acc_\ell(\check{x}, \check{y})} - g_i(p^\acc_\ell(\check{x}, \check{y})) - f_i(\hat{y}) + F_i(\hat{x}) \right] \\
           &+ \sum_{i = 1}^{m} \lambda_i(\check{x}, \check{y}) \left[ \innerp{\nabla f_i(\check{y})}{\check{y} - p^\acc_\ell(\hat{x}, \hat{y})} - g_i(p^\acc_\ell(\hat{x}, \hat{y})) - f_i(\check{y}) + F_i(\check{x}) \right]  
    .\end{align}
    From the definition~\cref{eq: p theta acc} of~$p^\acc_\ell$ and~$\theta^\acc_\ell$ and~\cref{eq: lambda}, we have
    \begin{align}
        \MoveEqLeft - \ell \norm*{p^\acc_\ell(\hat{x}, \hat{y}) - p^\acc_\ell(\check{x}, \check{y})}^2 \\
        \ge{}& \sum_{i = 1}^{m} \lambda_i(\check{x}, \check{y}) \left[ \innerp{\nabla f_i(\hat{y})}{p^\acc_\ell(\hat{x}, \hat{y}) - \hat{y}} + g_i(p^\acc_\ell(\hat{x}, \hat{y})) + f_i(\hat{y}) - F_i(\hat{x}) \right] \\ 
             &+ \sum_{i = 1}^{m} \lambda_i(\hat{x}, \hat{y}) \left[ \innerp{\nabla f_i(\check{y})}{p^\acc_\ell(\check{x}, \check{y}) - \check{y}} + g_i(p^\acc_\ell(\check{x}, \check{y})) + f_i(\check{y}) - F_i(\check{x}) \right] \\
             &+ \sum_{i = 1}^{m} \lambda_i(\hat{x}, \hat{y}) \left[ \innerp{\nabla f_i(\hat{y})}{\hat{y} - p^\acc_\ell(\check{x}, \check{y})} - g_i(p^\acc_\ell(\check{x}, \check{y})) - f_i(\hat{y}) + F_i(\hat{x}) \right] \\
             &+ \sum_{i = 1}^{m} \lambda_i(\check{x}, \check{y}) \left[ \innerp{\nabla f_i(\check{y})}{\check{y} - p^\acc_\ell(\hat{x}, \hat{y})} - g_i(p^\acc_\ell(\hat{x}, \hat{y})) - f_i(\check{y}) + F_i(\check{x}) \right] \\
             &- \frac{\ell}{2} 
             \begin{multlined}[t]
                 \Big[  \norm{p^\acc_\ell(\hat{x}, \hat{y}) - \check{y}}^2 - \norm{p^\acc_\ell(\hat{x}, \hat{y}) - \hat{y}}^2 \\
                  + \norm{p^\acc_\ell(\check{x}, \check{y}) - \hat{y}}^2 - \norm{p^\acc_\ell(\check{x}, \check{y}) - \check{y}}^2 \Big]
             \end{multlined} \\
        ={} & \sum_{i = 1}^{m} \lambda_i(\hat{x}, \hat{y})
        \begin{multlined}[t]
            [ \innerp{\nabla f_i(\hat{y})}{\hat{y}- \check{y}} + \innerp{\nabla f_i(\hat{y}) - \nabla f_i(\check{y})}{\check{y} - p^\acc_\ell(\check{x}, \check{y})} \\
            - f_i(\hat{y}) + f_i(\check{y}) + F_i(\hat{x}) - F_i(\check{x}) ] 
        \end{multlined} \\
        &+ \sum_{i = 1}^{m} \lambda_i(\check{x}, \check{y})
        \begin{multlined}[t]
            [ \innerp{\nabla f_i(\check{y})}{\check{y}- \hat{y}} + \innerp{\nabla f_i(\check{y}) - \nabla f_i(\hat{y})}{\hat{y} - p^\acc_\ell(\hat{x}, \hat{y})} \\
            - f_i(\check{y}) + f_i(\hat{y}) + F_i(\check{x}) - F_i(\hat{x}) ]
        \end{multlined} \\
        &- \ell \innerp{p^\acc_\ell(\hat{x}, \hat{y}) - p^\acc_\ell(\check{x}, \check{y})}{\hat{y} - \check{y}}
    .\end{align}
    Thus, \cref{eq: lambda} and Cauchy-Schwarz inequalities applied in each inner product that appears in the right-hand side of the above expression imply
    \begin{align}
        \MoveEqLeft - \ell \norm*{p^\acc_\ell(\hat{x}, \hat{y}) - p^\acc_\ell(\check{x}, \check{y})}^2 \\
        \ge{} & - 2 \max_{i = 1, \dots, m} \norm{\nabla f_i(\hat{y})} \norm{\hat{y} - \check{y}} \\
              &- \Big[ \norm{\hat{y} - p^\acc_\ell(\hat{x}, \hat{y})} + \norm{\check{y} - p^\acc_\ell(\check{x}, \check{y})} \Big] \max_{i = 1, \dots, m} \norm{\nabla f_i(\hat{y}) - \nabla f_i(\check{y})} \\
              &- 2 \max_{i = 1, \dots, m} \abs{f_i(\hat{y}) - f_i(\check{y})} - 2 \max_{i = 1, \dots, m} \abs{F_i(\hat{x}) - F_i(\check{x})} \\
              &- \ell \norm{p^\acc_\ell(\hat{x}, \hat{y}) - p^\acc_\ell(\check{x}, \check{y})} \norm{\hat{y} - \check{y}}
    .\end{align}
    Let us now show that each term of the right-hand side of the above inequality is bounded by a positive constant multiple of~$- \norm{\hat{x} - \check{x}}$ or~$- \norm{\hat{y} - \check{y}}$.
    The first term is direct because the boundedness of~$W$ implies~$\max_{i = 1, \dots, m} \norm{\nabla f_i(\hat{y})} < + \infty$.
    Since~$W$ is bounded and the objective function of~\cref{eq: acc prox subprob} is strongly convex,~$p^\acc_\ell(x, y)$ also belongs to some bounded set for all~$x, y \in W$, thus~$\norm{\hat{y} - p^\acc_\ell(\hat{x}, \hat{y})} < + \infty$ and~$\norm{\check{y} - p^\acc_\ell(\check{x}, \check{y})} < + \infty$.
    Thus, the Lipschitz continuity of~$\nabla f_i$ shows such a boundedness of the second term.
    Moreover, the locally Lipschitz continuity of~$f_i$ and~$F_i$ derived by the continuous differentiability of~$f_i$ and convexity~$F_i$ lead to the similar property for the third and fourth terms.
    Hence, $p^\acc_\ell$ is H\"{o}lder continuous with exponent~$1 / 2$ on~$W$.

    On the other hand, the definition~\cref{eq: p theta acc} of~$p^\acc_\ell$ and~$\theta^\acc_\ell$ gives
    \begin{align}
        \MoveEqLeft \theta^\acc_\ell(\hat{x}, \hat{y}) - \theta^\acc_\ell(\check{x}, \check{y}) \le \varphi^\acc_\ell(p^\acc_\ell(\check{x}, \check{y}); \hat{x}, \hat{y}) - \varphi^\acc_\ell(p^\acc_\ell(\check{x}, \check{y}); \check{x}, \check{y})  \\
        ={} & \max_{i = 1, \dots, m} [\innerp{\nabla f_i(\hat{y})}{p^\acc_\ell(\check{x}, \check{y}) - \hat{y}} + g_i(p^\acc_\ell(\check{x}, \check{y})) + f_i(\hat{y}) - F_i(\hat{x})] \\
            &- \max_{i = 1, \dots, m} [\innerp{\nabla f_i(\check{y})}{p^\acc_\ell(\check{x}, \check{y}) - \check{y}} + g_i(p^\acc_\ell(\check{x}, \check{y})) + f_i(\check{y}) - F_i(\check{x})] \\
           &+ \frac{\ell}{2} \Big[ \norm{p^\acc_\ell(\check{x}, \check{y}) - \hat{y}}^2 - \norm{p^\acc_\ell(\check{x}, \check{y}) - \check{y}}^2 \Big] \\
       \le{} & \max_{i = 1, \dots, m} 
       \begin{multlined}[t]
           [ \innerp{\nabla f_i(\check{y})}{\check{y} - \hat{y}} + \innerp{\nabla f_i(\hat{y}) - \nabla f_i(\check{y})}{p^\acc_\ell(\check{x}, \check{y}) - \hat{y}} \\
           + f_i(\hat{y}) - f_i(\check{y}) - F_i(\hat{x}) + F_i(\check{x}) ]
       \end{multlined} \\
             &+ \frac{\ell}{2} \innerp{2 p^\acc_\ell(\check{x}, \check{y}) - \hat{y} - \check{y}}{\check{y} - \hat{y}} \\
       \le{} & \max_{i = 1, \dots, m} \norm{\nabla f_i(\check{y})} \norm{\hat{y} - \check{y}} + \norm{\hat{y} - p^\acc_\ell(\check{x}, \check{y})} \max_{i = 1, \dots, m} \norm{\nabla f_i(\hat{y}) - \nabla f_i(\check{y})} \\
             &+ \max_{i = 1, \dots, m} \abs{f_i(\hat{y}) - f_i(\check{y})} + \max_{i = 1, \dots, m} \abs{F_i(\hat{x}) - F_i(\check{x})} \\
             &+ \frac{\ell}{2} \norm{2 p^\acc_\ell(\check{x}, \check{y}) - \hat{y} - \check{y}} \norm{\hat{y} - \check{y}}
   ,\end{align}
   where the second inequality follows from the relation~$\max_{i = 1, \dots, m} a_i - \max_{i = 1, \dots, m} b_i \le \max_{i = 1, \dots, m} (a_i - b_i)$ for all~$a, b \in \setR^m$, and the third inequality comes from~\cref{eq: lambda} and Cauchy-Schwarz inequalites.
   Since the above inequality holds even if we interchange~$(\hat{x}, \hat{y})$ and~$(\check{x}, \check{y})$, we can show the Lipschitz continuity of~$\theta^\acc_\ell$ on~$W$ in the same way as in the previous paragraph.
     \end{proof}
     Note that the H\"{o}lder exponent~$1 / 2$ mentioned in \cref{enum: continuity} is optimal, i.e., for some~$F_i$,~$p^\acc_\ell$ is not H\"{o}lder continuous with exponent~$\alpha > 1 / 2$.
     In fact, this result was also proved for multiobjective steepest direction in~\cite{Svaiter2018}.

     \Cref{prop: acc prox termination} suggests that we can use~$\norm*{p^\acc_\ell(x, y) - y}_\infty < \varepsilon$ for some~$\varepsilon > 0$ as a stopping criteria.
Now, we state below the proposed algorithm.

\begin{algorithm}[hbtp]
    \caption{Accelerated proximal gradient method for multiobjective optimization}
    \label{alg: acc-pgm}
    \begin{algorithmic}[1]
        \Require Set~$x^0 = y^1 \in \dom F, \ell \ge L, \varepsilon > 0$.
        \Ensure $x^\ast$: A weakly Pareto optimal point
        \State $k \gets 1$
        \State $t_1 \gets 1$ \label{line: t ini}
        \Loop
        \State $x^k \gets p^\acc_\ell(x^{k - 1}, y^k)$ \label{line: x}, where $p^\acc_\ell$ is defined by~\cref{eq: p theta acc}
        \If{$\norm*{x^k - y^k}_\infty < \varepsilon$} \label{line: stop}
        \State \Return $x^k$
        \EndIf
        \State $t_{k + 1} \gets \sqrt{t_k^2 + 1 / 4} + 1/2$ \label{line: t rr}
        \State $\gamma_k \gets (t_k - 1) / t_{k + 1}$ \label{line: gamma}
        \State $y^{k + 1} \gets x^k + \gamma_k (x^k - x^{k - 1})$ \label{line: y}
        \State $k \gets k + 1$
        \EndLoop
    \end{algorithmic}
\end{algorithm}

We show below some properties of~$\setof*{t_k}$ and~$\setof*{\gamma_k}$, related to stepsizes.
\begin{lemma} \label{lem: t}
    Let~$\setof*{t_k}$ and~$\setof*{\gamma_k}$ be defined by \cref{line: t ini,line: t rr,line: gamma} in \cref{alg: acc-pgm}.
    Then, the following inequalities hold for all~$k \ge 1$: \\
    \begin{enumerate*}
        \item $t_{k + 1} \ge t_k + 1 / 2$ and~$t_k \ge (k + 1) / 2$; \label{enum: t geq}
        \item $t_k^2 - t_{k + 1}^2 + t_{k + 1} = 0$; \label{enum: t eq}
        \item $1 - \gamma_k^2 \ge \dfrac{1}{t_k}$. \label{enum: t moment}
    \end{enumerate*}
\end{lemma}
\begin{proof}
    \subCref{enum: t geq}:
    From the definition of~$\setof*{t_k}$, we have
    \[ \label{eq: t cs} 
        t_{k + 1} = \sqrt{t_k^2 + \frac{1}{4}} + \frac{1}{2}
        \ge t_k + \frac{1}{2}
    .\]
    Applying the above inequality recursively, we obtain
    \[
        t_k \ge t_1 + \frac{k - 1}{2} = \frac{k + 1}{2}
    .\]

    \subCref{enum: t eq}:
    An easy computation shows that
    \begin{align}
            t_k^2 - t_{k + 1}^2 + t_{k + 1} &= t_k^2 - \left[ \sqrt{t_k^2 + \frac{1}{4}} + \frac{1}{2} \right]^2 + \sqrt{t_k^2 + \frac{1}{4}} + \frac{1}{2} = 0
    .\end{align}

    \subCref{enum: t moment}:
    \subCref{enum: t geq} of this lemma implies that~$t_{k + 1} > t_k \ge 1$.
    Thus, the definition of~$\gamma_k$ leads to
    \[
        1 - \gamma_k^2 = 1 - \left( \frac{t_k - 1}{t_{k + 1}} \right)^2 \ge 1 - \left( \frac{t_k - 1}{t_k} \right)^2
        = \frac{2 t_k - 1}{t_k^2} \ge \frac{2 t_k - t_k}{t_k^2} = \frac{1}{t_k}
    .\]
\end{proof}

We end this section by noting some remarks about the proposed algorithm.
\begin{remark} \label{rem: acc-pgm}
    \begin{enumerate}
        \item When~$m = 1$, we can remove the term~$f_i(y) - F_i(x)$ from the subproblem~\cref{eq: acc prox subprob}, so \cref{alg: acc-pgm} corresponds to the Fast Iterative Shrinkage-Thresholding Algorithm (FISTA)~\cite{Beck2009} for single-objective optimization.
        \item \Cref{alg: acc-pgm} produces two sequences~$x^k$ and~$y^k$, in a similar way to the single-objective FISTA. In particular, the stopping condition (Step~\ref{line: stop}), the momentum update (Steps~\ref{line: t rr} and~\ref{line: gamma}), and the update of the iterate (Step~\ref{line: y}) are actually equivalent to the single-objective case.
        \item Since~$x \in \dom F$ implies~$p^\acc_\ell(x, y) \in \dom F$, every~$x^k$ computed by the above algorithm is in~$\dom F$.
            However,~$y^k$ is not necessarily in~$\dom F$.
        \item Since $y^1 = x^0$, it follows from~\cref{eq: acc prox subprob optimality} that~$\theta^\acc_\ell(x^0, y^1) \le 0$, but the inequality~$\theta^\acc_\ell(x^{k - 1}, y^k) \le 0$ does not necessarily hold for~$k \ge 2$.
        \item Like \cref{enum: other methods}, \cref{alg: acc-pgm} induces the accelerated versions of first-order algorithms such as the steepest descent~\cite{Fliege2000}, proximal point~\cite{Bonnel2005}, and projected gradient methods~\cite{Grana-Drummond2004}.
        \item Like \cref{enum: backtracking}, even if it is difficult to estimate~$L$, we can update the constant~$\ell$ to satisfy~$F_i(p^\acc_\ell(x^{k - 1}, y^k)) - F_i(x^{k - 1}) \le \theta^\acc_\ell(x^{k - 1}, y^k)$ for all~$i = 1, \dots, m$ in each iteration by a finite number of backtracking steps.
            Moreover, we can restrict the assumption of~$\nabla f_i$'s Lipschitz continuity on the level set~$\level_F(F(x^0))$ without affecting the analysis in the subsequent sections.
    \end{enumerate}
\end{remark}

\ifmain
\else
\fi

\section{Convergence rate} \label{sec: convergence rate}
This section shows that \cref{alg: acc-pgm} has a convergence rate of~$O(1/k^2)$ under the same assumptions used in the complexity analysis of \cref{alg: pgm}.
As it is expected, this rate is better than the one obtained for \cref{alg: pgm}.

Let us first define some functions below, that will be useful for our complexity analysis.
For~$k \ge 0$, let~$\sigma_k \colon \setR^n \to \setR \cup \setof*{- \infty}$ and~$\rho_k \colon \setR^n \to \setR$ be defined by
\[ \label{eq: sigma rho}
\begin{gathered} 
    \sigma_k(z) \coloneqq \min_{i = 1, \dots, m}\left[ F_i(x^k) - F_i(z) \right], \\
        \rho_k(z) \coloneqq \norm*{t_{k + 1} x^{k + 1} - (t_{k + 1} - 1) x^k - z}^2  
,\end{gathered}
\]
respectively.
We present a lemma on~$\sigma_k$ that will be helpful in the subsequent discussions.
\begin{lemma} \label{lem: sigma}
    Let~$\setof{x^k}$ and~$\setof{y^k}$ be sequences generated by \cref{alg: acc-pgm}.
    Then, the following inequalities hold for all~$z \in \setR^n$ and~$k \ge 0$:
    \begin{gather}
        \begin{multlined}
            \sigma_{k + 1}(z) \le - \frac{\ell}{2} \left[ 2 \innerp*{x^{k + 1} - y^{k + 1}}{y^{k + 1} - z} + \norm*{x^{k + 1} - y^{k + 1}}^2 \right] \\
            - \frac{\ell - L}{2} \norm*{x^{k + 1} - y^{k + 1}}^2
        ,\end{multlined} \label{eq: lemsigma1} \\
        \begin{multlined}
            \sigma_k(z) - \sigma_{k + 1}(z) \ge \frac{\ell}{2} \left[ 2 \innerp*{x^{k + 1} - y^{k + 1}}{y^{k + 1} - x^k} + \norm*{x^{k + 1} - y^{k + 1}}^2 \right] \\
            + \frac{\ell - L}{2} \norm*{x^{k + 1} - y^{k + 1}}^2
        .\end{multlined} \label{eq: lemsigma2}
    \end{gather}
\end{lemma}
\begin{proof}
    Suppose that~$z \in \setR^n$ and~$k \ge 0$.
    Recall that there exist~$\eta(x^k, y^{k + 1}) \in \partial g(x^{k + 1})$ and a Lagrange multiplier~$\lambda(x^k, y^{k + 1}) \in \setR^m$ that satisfy the KKT condition~\cref{eq: kkt} for the subproblem~\cref{eq: acc prox subprob}.
    From the definition~\cref{eq: sigma rho} of~$\sigma_{k + 1}$, we get
    \[
        \sigma_{k + 1}(z) = \min_{i = 1, \dots, m} \left[ F_i(x^{k + 1}) - F_i(z) \right]
        \le \sum_{i = 1}^m \lambda_i(x^k, y^{k + 1}) \left[ F_i(x^{k + 1}) - F_i(z) \right]
    .\]
    where the inequality follows from~\cref{eq: lambda}.
    Taking~$p = x^{k + 1}, q = y^{k + 1}$, and~$r = z$ in~\cref{eq: descent mid}, we have
    \begin{multline}
        \sigma_{k + 1}(z)
    \le \sum_{i = 1}^m \lambda_i(x^k, y^{k + 1}) \left[ \innerp*{\nabla f_i(y^{k + 1})}{x^{k + 1} - y^{k + 1}} + g_i(x^{k + 1}) \right.\\ \left. + f_i(y^{k + 1}) - F_i(z) \right] 
        + \frac{L}{2} \norm*{x^{k + 1} - y^{k + 1}}^2
    .\end{multline}
    Hence, the convexity of~$f_i$ and~$g_i$ yields
    \begin{align}
            &\sigma_{k + 1}(z)\\
            &\begin{multlined}
                \le \sum_{i = 1}^m \lambda_i(x^k, y^{k + 1}) \left[ \innerp*{\nabla f_i(y^{k + 1})}{x^{k + 1} - y^{k + 1}} + \innerp*{\nabla f_i(y^{k + 1})}{y^{k + 1} - z} \right. \\
                + \left. \innerp*{\eta_i(x^k, y^{k + 1})}{x^{k + 1} - z} \right] + \frac{L}{2} \norm*{x^{k + 1} - y^{k + 1}}^2
            \end{multlined} \\
            &= \sum_{i = 1}^m \lambda_i(x^k, y^{k + 1}) \innerp*{\nabla f_i(y^{k + 1}) + \eta_i(x^k, y^{k + 1})}{x^{k + 1} - z} + \frac{L}{2} \norm*{x^{k + 1} - y^{k + 1}}^2
    .\end{align}
    Using~\cref{eq: optimal} with~$x = x^k$ and~$y = y^{k + 1}$ and from the fact that~$x^{k + 1} = p^\acc_\ell(x^k, y^{k + 1})$ (see \cref{line: x} of \cref{alg: acc-pgm}), we obtain
    \[
        \sigma_{k + 1}(z) \le - \ell \innerp*{x^{k + 1} - y^{k + 1}}{x^{k + 1} - z} + \frac{L}{2} \norm*{x^{k + 1} - y^{k + 1}}^2
    .\]
    Moreover, simple calculations show that
    \begin{align}
            &\sigma_{k + 1}(z) \\
            &\le - \frac{\ell}{2} \left[ 2 \innerp*{x^{k + 1} - y^{k + 1}}{x^{k + 1} - z} - \norm*{x^{k + 1} - y^{k + 1}}^2 \right] - \frac{\ell - L}{2} \norm*{x^{k + 1} - y^{k + 1}}^2 \\
            &= - \frac{\ell}{2} \left[ 2 \innerp*{x^{k + 1} - y^{k + 1}}{y^{k + 1} - z} + \norm*{x^{k + 1} - y^{k + 1}}^2 \right] - \frac{\ell - L}{2} \norm*{x^{k + 1} - y^{k + 1}}^2
    ,\end{align}
    which completes the proof of~\cref{eq: lemsigma1}.

    Now, let us show inequality~\cref{eq: lemsigma2}.
    Again from the definition~\cref{eq: sigma rho} of~$\sigma_k$, we obtain
    \[ \label{eq: lemsigma tmp}
        \begin{aligned}
            \MoveEqLeft \sigma_k(z) - \sigma_{k + 1}(z) \\
        &= \min_{i = 1, \dots, m}\left[ F_i(x^k) - F_i(z) \right] - \min_{i = 1, \dots, m}\left[ F_i(x^{k + 1}) - F_i(z) \right] \\
        &\ge - \max_{i = 1, \dots, m}\left[ F_i(x^{k + 1}) - F_i(x^k) \right]
        ,\end{aligned}
    \] 
    where the inequality holds because
    \[
        \min_{i = 1, \dots, m} \left( u_i + v_i \right) - \min_{i = 1, \dots, m} u_i \ge \min_{i = 1, \dots, m} v_i \forallcondition{u, v \in \setR^m}
    .\]
    Letting~$p = x^{k + 1}, q = y^{k + 1}$, and~$r = x^k$ in~\cref{eq: descent mid}, we have
    \begin{align}
            \MoveEqLeft
            \sigma_k(z) - \sigma_{k + 1}(z) \\
            &\begin{multlined}
                \ge - \max_{i = 1, \dots, m} \left[ \innerp*{\nabla f_i(y^{k + 1})}{x^{k + 1} - y^{k + 1}} + g_i(x^{k + 1}) + f_i(y^{k + 1}) \right. \\
                - \left. F_i(x^k) \right] - \frac{L}{2} \norm*{x^{k + 1} - y^{k + 1}}^2
            \end{multlined} \\
            &\begin{multlined}
                = - \sum_{i = 1}^m \lambda_i(x^k, y^{k + 1}) \left[ \innerp*{\nabla f_i(y^{k + 1})}{x^{k + 1} - y^{k + 1}} + g_i(x^{k + 1}) \right. \\
                + \left. f_i(y^{k + 1}) - F_i(x^k) \right] - \frac{L}{2} \norm*{x^{k + 1} - y^{k + 1}}^2 
            \end{multlined} \\
            &\begin{multlined}
                = - \sum_{i = 1}^m \lambda_i(x^k, y^{k + 1}) \left[ \innerp*{\nabla f_i(y^{k + 1})}{x^k - y^{k + 1}} + f_i(y^{k + 1}) - f_i(x^k) \right] \\
                - \sum_{i = 1}^m \lambda_i(x^k, y^{k + 1}) \left[ \innerp*{\nabla f_i(y^{k + 1})}{x^{k + 1} - x^k} + g_i(x^{k + 1}) - g_i(x^k) \right] \\
                - \frac{L}{2} \norm*{x^{k + 1} - y^{k + 1}}^2
            ,\end{multlined}
    \end{align}
    where the first equality comes from~\cref{eq: lambda}, and the second one follows by taking~$x^{k + 1} - y^{k + 1} = (x^k - y^{k + 1}) + (x^{k + 1} - x^k)$.
    From the convexity of~$f_i$, the first term of the above expression is nonnegative.
    Moreover, the convexity of~$g_i$ shows that
    \begin{multline}
        \sigma_k(z) - \sigma_{k + 1}(z) \\
        \ge - \sum_{i = 1}^m {\lambda_i(x^k, y^{k + 1}) \innerp*{\nabla f_i(y^{k + 1}) + \eta_i(x^k, y^{k + 1})}{x^{k + 1} - x^k}} - \frac{L}{2} \norm*{x^{k + 1} - y^{k + 1}}^2
    .\end{multline}
    Thus,~\cref{eq: optimal} with~$(x, y) = (x^k, y^{k + 1})$ and direct calculations prove that
    \begin{align}
        &\sigma_k(z) - \sigma_{k + 1}(z) \\
        &\ge \ell \innerp*{x^{k + 1} - y^{k + 1}}{x^{k + 1} - x^k} - \frac{L}{2} \norm*{x^{k + 1} - y^{k + 1}}^2 \\
        &= \frac{\ell}{2} \left[ 2 \innerp*{x^{k + 1} - y^{k + 1}}{x^{k + 1} - x^k} - \norm*{x^{k + 1} - y^{k + 1}}^2 \right] + \frac{\ell - L}{2} \norm*{x^{k + 1} - y^{k + 1}}^2 \\
        &= \frac{\ell}{2} \left[ 2 \innerp*{x^{k + 1} - y^{k + 1}}{y^{k + 1} - x^k} + \norm*{x^{k + 1} - y^{k + 1}}^2 \right] + \frac{\ell - L}{2} \norm*{x^{k + 1} - y^{k + 1}}^2
    .\end{align}
\end{proof}

We can also show the following corollary of \cref{lem: sigma}~\cref{eq: lemsigma2}.
\begin{corollary} \label{cor: sigma}
    Let~$\setof{x^k}$ and~$\setof{y^k}$ be sequences generated by \cref{alg: acc-pgm}.
    Then, we have
    \begin{multline}
        \sigma_{k_1}(z) - \sigma_{k_2}(z) \\
        \ge \frac{\ell}{2} \left[ \norm*{x^{k_2} - x^{k_2 - 1}}^2 - \norm*{x^{k_1} - x^{k_1 - 1}}^2 + \sum_{k = k_1}^{k_2 - 1} \frac{1}{t_k} \norm*{x^k - x^{k - 1}}^2 \right]
    \end{multline}
    for any~$k_2 \ge k_1 \ge 1$.
\end{corollary}
\begin{proof}
    Let~$k \ge 1$.
    Since~$\ell \ge L$, \Cref{lem: sigma}~\cref{eq: lemsigma2} implies
    \begin{align}
            \sigma_k(z) - \sigma_{k + 1}(z) &\ge \frac{\ell}{2} \left[ 2 \innerp*{x^{k + 1} - y^{k + 1}}{y^{k + 1} - x^k} + \norm*{x^{k + 1} - y^{k + 1}}^2 \right] \\
        &= \frac{\ell}{2} \left[ \norm*{x^{k + 1} - x^k}^2 - \norm*{y^{k + 1} - x^k}^2 \right]
    ,\end{align}
    where the equality holds from~\cref{eq: Pythagoras} with~$(a, b, c) = (y^{k + 1}, x^{k + 1}, x^k)$.
    Hence, the definition of~$y^{k + 1}$ given in \cref{line: y} of \cref{alg: acc-pgm} yields
    \[
        \sigma_k(z) - \sigma_{k + 1}(z) \ge \frac{\ell}{2} \left[ \norm*{x^{k + 1} - x^k}^2 - \gamma_k^2 \norm*{x^k - x^{k - 1}}^2 \right]
    .\]
    Applying this inequality repeatedly, we have
    \begin{multline}
        \sigma_{k_1}(z) - \sigma_{k_2}(z) \\
        \ge \frac{\ell}{2} \left[ \norm*{x^{k_2} - x^{k_2 - 1}}^2 - \norm*{x^{k_1} - x^{k_1 - 1}}^2 + \sum_{k = k_1}^{k_2 - 1} \left( 1 - \gamma_k^2 \right)  \norm*{x^k - x^{k - 1}}^2 \right] 
    .\end{multline}
    Using \cref{enum: t moment}, we get the desired inequality.
\end{proof}

Before analyzing the convergence rate of \cref{alg: acc-pgm}, we show that the objective function values at~$x^k$ for any~$k \ge 0$ never exceed the ones at the initial point, that is,~$\setof{x^k}$ belongs to the level set~$\level_F(F(x^0))$ (see~\cref{eq: level set} for the definition of~$\level_F$).
However, note that \cref{alg: acc-pgm} does not guarantee the monotonically decreasing property~$F(x^{k + 1}) \le F(x^k)$.
\begin{theorem} \label{thm: leq ini}
    \cref{alg: acc-pgm} generates a sequence~$\setof{x^k}$ such that
    \[
        F_i(x^k) \le F_i(x^0) \forallcondition{i = 1, \dots, m, k \ge 0}
    .\]
\end{theorem}
\begin{proof}
    Let~$i = 1, \dots, m$ and~$p \ge 1$.
    Then, we have
    \[
        F_i(x^p) - F_i(x^{p + 1}) \ge - \max_{i = 1, \dots, m} \left[ F_i(x^{p + 1}) - F_i(x^p) \right]
    .\]
    Noting that the right-hand side of the above inequality is the same as~\cref{eq: lemsigma tmp}, with similar arguments used in the proof of~\cref{eq: lemsigma2} in \cref{lem: sigma}, we obtain
    \begin{multline} \label{eq: F - F}
        F_i(x^p) - F_i(x^{p + 1})
        \ge \frac{\ell}{2} \left[ 2 \innerp*{x^{p + 1} - y^{p + 1}}{y^{p + 1} - x^p} + \norm*{x^{p + 1} - y^{p + 1}}^2 \right] \\ + \frac{\ell - L}{2} \norm*{x^{p + 1} - y^{p + 1}}^2
    .\end{multline}
    Note that this inequality also holds for~$p = 0$.
    Again, in the same way as in the proof of \cref{cor: sigma}, we get
    \[
        F_i(x^1) - F_i(x^k) \ge \frac{\ell}{2} \left[ \norm*{x^k - x^{k - 1}}^2 - \norm*{x^1 - x^0}^2 + \sum_{p = 1}^{k - 1} \frac{1}{t_p} \norm*{x^p - x^{p - 1}}^2 \right]
    .\]
    Since~$t_1 = 1$, the above inequality reduces to
    \[
        F_i(x^1) - F_i(x^k) \ge \frac{\ell}{2} \left[ \norm*{x^k - x^{k - 1}}^2 + \sum_{p = 2}^{k - 1} \frac{1}{t_p} \norm*{x^p - x^{p - 1}}^2 \right] \ge 0
    .\] 
    Moreover,~\cref{eq: F - F} with~$p = 0$ and the fact that~$y^1 = x^0$ imply~$F_i(x^1) \le F_i(x^0)$, so we can conclude that~$F_i(x^k) \le F_i(x^0)$.
\end{proof}

The following result provides the fundamental relation for our convergence rate analysis.
\begin{lemma} \label{lem: key relation}
    Let~$\setof{x^k}$ and~$\setof{y^k}$ be sequences generated by \cref{alg: acc-pgm}.
    Also, let~$\sigma_k$ and~$\rho_k$ be defined by~\cref{eq: sigma rho}.
    Then, we have
    \[
        t_{k + 1}^2 \sigma_{k + 1}(z) + \frac{\ell}{2} \rho_k(z) + \frac{\ell - L}{2} \sum_{p = 1}^{k} t_{p + 1}^2 \norm*{x^{p + 1} - y^{p + 1}}^2 \le 
        \frac{\ell}{2} \norm*{x^0 - z}^2 
    \]
    for all~$k \ge 0$ and~$z \in \setR^n$.
\end{lemma}
\begin{proof}
    Let~$p \ge 1$ and~$z \in \setR^n$.
    Recall from \cref{lem: sigma} that
    \begin{gather}
        \begin{multlined}
        - \sigma_{p + 1}(z) \ge \frac{\ell}{2} \left[ 2 \innerp*{x^{p + 1} - y^{p + 1}}{y^{p + 1} - z} + \norm*{x^{p + 1} - y^{p + 1}}^2 \right] \\
        + \frac{\ell - L}{2} \norm*{x^{p + 1} - y^{p + 1}}^2,
        \end{multlined} \\
        \begin{multlined}
            \sigma_p(z) - \sigma_{p + 1}(z) \ge \frac{\ell}{2} \left[ 2 \innerp*{x^{p + 1} - y^{p + 1}}{y^{p + 1} - x^p} + \norm*{x^{p + 1} - y^{p + 1}}^2 \right] \\
            + \frac{\ell - L}{2} \norm*{x^{p + 1} - y^{p + 1}}^2
        .\end{multlined}
        \end{gather}
        To get a relation between $\sigma_p(z)$ and $\sigma_{p + 1}(z)$, we multiply the second inequality above by $(t_{p + 1} - 1)$ and add it to the first one:
    \begin{multline}
        (t_{p + 1} - 1) \sigma_p(z) - t_{p + 1} \sigma_{p + 1}(z) \\
        \ge \frac{\ell}{2} \left[ t_{p + 1} \norm*{x^{p + 1} - y^{p + 1}}^2 + 2 \innerp*{x^{p + 1} - y^{p + 1}}{t_{p + 1} y^{p + 1} - (t_{p + 1} - 1)x^p - z} \right] \\
        + \frac{\ell - L}{2} t_{p + 1} \norm*{x^{p + 1} - y^{p + 1}}^2
    .\end{multline}
    Multiplying this inequality by~$t_{p + 1}$ and using the relation~$t_p^2 = t_{p + 1}^2 - t_{p + 1}$ (cf. \cref{enum: t eq}), we get
    \begin{multline}
        t_p^2 \sigma_p(z) - t_{p + 1}^2 \sigma_{p + 1}(z) \\
        \ge \frac{\ell}{2} \left[ \norm*{t_{p + 1} (x^{p + 1} - y^{p + 1})}^2
        + 2 t_{p + 1} \innerp*{x^{p + 1} - y^{p + 1}}{t_{p + 1} y^{p + 1} - (t_{p + 1} - 1)x^p - z} \right] \\
        + \frac{\ell - L}{2} t_{p + 1}^2 \norm*{x^{p + 1} - y^{p + 1}}^2
    .\end{multline}
    Applying~\cref{eq: Pythagoras} with~$(a, b, c) = \left(t_{p + 1} y^{p + 1}, t_{p + 1} x^{p + 1}, (t_{p + 1} - 1) x^p + z \right)$ to the right-hand side of the last inequality, we get
    \begin{multline}
        t_p^2 \sigma_{p + 1}(z) - t_{p + 1}^2 \sigma_p(z) \\
        \ge \frac{\ell}{2} \left[ \norm*{t_{p + 1} x^{p + 1} - (t_{p + 1} - 1) x^p - z}^2 - \norm*{t_{p + 1} y^{p + 1} - (t_{p + 1} - 1) x^p - z}^2 \right] \\
        + \frac{\ell - L}{2} t_{p + 1}^2 \norm*{x^{p + 1} - y^{p + 1}}^2
    .\end{multline}
    Recall that~$\rho_p(z) \coloneqq \norm*{t_{p + 1} x^{p + 1} - (t_{p + 1} - 1) x^p - z}^2$.
    Then, from the definition of~$y^p$ defined in \cref{line: y} of \cref{alg: acc-pgm}, we get
    \[
        t_p^2 \sigma_p(z) - t_{p + 1}^2 \sigma_{p + 1}(z)
        \ge \frac{\ell}{2} \left[ \rho_p(z) - \rho_{p - 1}(z) \right] + \frac{\ell - L}{2} t_{p + 1}^2 \norm*{x^{p + 1} - y^{p + 1}}^2
    .\]
    Now, let~$k \ge 0$.
    Adding the above inequality from~$p = 0$ to~$p = k$ and using~$t_1 = 1$ and~$\rho_0(z) = \norm*{x^1 - z}^2$, we have
    \[ \label{eq: key relation tmp}
        \sigma_1(z) - t_{k + 1}^2 \sigma_{k + 1}(z)
        \ge \frac{\ell}{2} \left[ \rho_{k}(z) - \norm*{x^1 - z}^2 \right] + \frac{\ell - L}{2} \sum_{p = 1}^k t_{k + 1}^2 \norm*{x^{k + 1} - y^{k + 1}}^2
    .\]
    \Cref{lem: t}~\cref{eq: lemsigma1} with~$k = 0$ and~$y^1 = x^0$ lead to
    \begin{align}
        \sigma_1(z) &\le - \frac{\ell}{2} \left[ \norm*{x^1 - z}^2 - \norm*{x^0 - z}^2 \right] - \frac{\ell - L}{2} \norm*{x^1 - y^1}^2 \\
                    &\le - \frac{\ell}{2} \left[ \norm*{x^1 - z}^2 - \norm*{x^0 - z}^2 \right]
    ,\end{align}
    where the second inequality follows since~$\ell \ge L$.
    From the above two inequalities, we can derive the desired inequality.
\end{proof}

Finally, using \cref{lem: key relation}, we can evaluate the convergence rate of \cref{alg: acc-pgm} with the following theorem.
\begin{theorem} \label{thm: acc conv rate}
    Under \cref{asm: bound}, \cref{alg: acc-pgm} generates a sequence~$\setof{x^k}$ such that
    \[
        u_0(x^k) \le \frac{2 \ell R}{(k + 1)^2}
    ,\]
    where~$R \ge 0$ is given in~\cref{eq: R}, and~$u_0$ is a merit function defined by~\cref{eq: u_0}.
\end{theorem}
\begin{proof}    
    Let~$k \ge 0$.
    Since~$\rho_k(z) \ge 0$, \cref{lem: key relation} gives
    \[
        t_{k + 1}^2 \sigma_{k + 1}(z) \le \frac{\ell}{2} \norm*{x^0 - z}^2
    .\] 
    It follows from \cref{enum: t geq} that
    \[
        \sigma_{k + 1}(z) \le \frac{2 \ell}{(k + 2)^2} \norm*{x^0 - z}^2
    .\] 
    With similar arguments used in the proof of \cref{thm: conv rate} (see~\cite[Theorem 5.2]{Tanabe2023}), we get the desired inequality.
\end{proof}

We end this section by showing that the global convergence of \cref{alg: acc-pgm}, in terms of weak Pareto optimality, is also guaranteed by using the above compexity result.
\begin{corollary}
    Suppose that \cref{asm: bound} holds.
    Then, every accumulation point of the sequence~$\setof{x^k}$ generated by \cref{alg: acc-pgm} is weakly Pareto optimal for~\cref{eq: MOP}.
    In particular, if the level set~$\level_F(F(x^0))$ is bounded, then~$\setof{x^k}$ has accumulation points, and they are all weakly Pareto optimal.

    Moreover, if each~$F_i$ is strictly convex, then the accumulation points are Pareto optimum, i.e., there does not exist any points with the same or smaller objective function values and with at least one objective function value being strictly smaller.
\end{corollary}
\begin{proof}
    The first claim is clear from the lower-semicontinuity of~$F_i$ for all~$i = 1, \dots, m$ as well as \cref{thm: merit Pareto,thm: acc conv rate}, and the second one is easy since \cref{thm: leq ini} holds.
    The third is also obvious from the relationship between weak Pareto and Pareto optimalities~\cite[Lemma 2.2]{Tanabe2019}.
\end{proof}

\ifmain
\else
\fi

\section{Efficient computation of the subproblem via its dual} \label{sec: subproblem}
In the previous section, we proved global convergence and complexity results of \cref{alg: acc-pgm}.
Now, we want to show how practical is the proposed method.
In particular, we now discuss a way of computing the subproblem~\cref{eq: acc prox subprob}.
First, define
\[ \label{eq: psi}
    \psi_i(z; x, y) \coloneqq \innerp*{\nabla f_i(y)}{z - y} + g_i(z) + f_i(y) - F_i(x) + \frac{\ell}{2} \norm*{z - y}^2
\]
for all~$i = 1, \dots, m$.
Then, fixing some~$\ell \ge L$, we can rewrite the objective function~$\varphi^\acc_\ell(z; x, y)$ of~\cref{eq: acc prox subprob} as
\[
    \varphi^\acc_\ell(z; x, y) = \max_{i = 1, \dots, m} \psi_i(z; x, y)
.\] 
Recall that~$\Delta^m \subseteq \setR^m$ represents the standard simplex~\cref{eq: simplex}.
Since~$\max_{i = 1, \dots, m} q_i = \max_{\lambda \in \Delta^m} \sum_{i = 1}^{m} \lambda_i q_i$ for any~$q \in \setR^m$, we get
\[
    \varphi^\acc_\ell(z; x, y) = \max_{\lambda \in \Delta^m} \sum_{i = 1}^{m} \lambda_i \psi_i(z; x, y)
.\] 
Then, the subproblem~\cref{eq: acc prox subprob} reduces to the following minimax problem:
\[ \label{eq: min max}
    \min_{z \in \setR^n} \max_{\lambda \in \Delta^m} \quad \sum_{i = 1}^{m} \lambda_i \psi_i(z; x, y)
.\] 
We can see that~$\setR^n$ is convex,~$\Delta^m$ is compact and convex, and~$\sum_{i = 1}^{m} \lambda_i \psi_i(z; x, y)$ is convex for~$z$ and concave for~$\lambda$.
Therefore, Sion's minimax theorem~\cite{Sion1958} shows that the above problem is equivalent to
\[ \label{eq: max min}
    \max_{\lambda \in \Delta^m} \min_{z \in \setR^n} \quad \sum_{i = 1}^{m} \lambda_i \psi_i(z; x, y)
.\] 
The definition~\cref{eq: psi} of~$\psi_i$ yields
\begin{align}
        \min_{z \in \setR^n} \sum_{i = 1}^{m} \lambda_i \psi_i(z; x, y) ={}& \min_{z \in \setR^n} \left[ \sum_{i = 1}^{m} \lambda_i g_i(z) + \frac{\ell}{2} \norm*{z - y + \frac{1}{\ell} \sum_{i = 1}^{m} \lambda_i \nabla f_i(y)}^2 \right] \\
                                                                         &- \frac{1}{2 \ell} \norm*{\sum_{i = 1}^{m} \lambda_i \nabla f_i(y)}^2 + \sum_{i = 1}^{m} \lambda_i \left\{ f_i(y) - F_i(x) \right\} \\
        ={}& \ell \envelope_{\frac{1}{\ell}\sum\limits_{i = 1}^{m} \lambda_i g_i}\left( y - \frac{1}{\ell} \sum_{i = 1}^{m} \lambda_i \nabla f_i(y) \right) \\
                                                                         &- \frac{1}{2 \ell} \norm*{\sum_{i = 1}^{m} \lambda_i \nabla f_i(y)}^2 + \sum_{i = 1}^{m} \lambda_i \left\{ f_i(y) - F_i(x) \right\}
,\end{align}
where~$\envelope$ is the Moreau envelope~\cref{eq: Moreau envelope}.
Based on the discussion above, we obtain the dual problem of~\cref{eq: acc prox subprob} as follows:
\[ \label{eq: dual}
    \begin{aligned}
        \max_{\lambda \in \setR^m} &&& \omega(\lambda) \\
        \st &&& \lambda \ge 0 \eqand \sum_{i = 1}^{m} \lambda_i = 1 
    ,\end{aligned}
\] 
where
\[ \label{eq: omega}
    \begin{aligned}
        \omega(\lambda) \coloneqq{}& \ell \envelope_{\frac{1}{\ell} \sum\limits_{i = 1}^{m} \lambda_i g_i}\left( y - \frac{1}{\ell} \sum_{i = 1}^{m} \lambda_i \nabla f_i(y) \right) \\&- \frac{1}{2 \ell} \norm*{\sum_{i = 1}^{m} \lambda_i \nabla f_i(y)}^2 + \sum_{i = 1}^{m} \lambda_i \left[ f_i(y) - F_i(x) \right]
    .\end{aligned}
\] 
If we can find the global optimal solution~$\lambda^\ast$ of this dual problem~\cref{eq: dual}, we can construct the optimal solution~$z^\ast$ of the original subproblem~\cref{eq: acc prox subprob} as
\[
    z^\ast = \prox_{\frac{1}{\ell} \sum\limits_{i = 1}^{m} \lambda_i^\ast g_i}\left( y - \frac{1}{\ell}\sum_{i = 1}^{m} \lambda_i^\ast \nabla f_i(y) \right) 
,\] 
where~$\prox$ denotes the proximal operator~\cref{eq: proximal operator}.
This is because the equivalence between~\cref{eq: min max,eq: max min} induces
\[
\sum_{i = 1}^{m} \lambda^\ast_i \psi_i(z^\ast; x, y) = \max_{\lambda \in \Delta^n} \min_{z \in \setR^n} \sum_{i = 1}^{m} \lambda_i \psi_i(z; x, y) = \min_{z \in \setR^n} \max_{\lambda \in \Delta^n} \sum_{i = 1}^{m} \lambda_i \psi_i(z; x, y)
,\]
which means that~$z^\ast$ attains the minimum in~\cref{eq: min max}.
Since~$\sum_{i = 1}^{m} \lambda_i \psi_i(z; x, y)$ is concave for~$\lambda$, it is clear that~$\omega(\lambda) = \min_{z \in \setR^n}\lambda_i \psi_i(z; x, y)$ is concave.
Furthermore,~$\omega$ is differentiable, as the following theorem shows.
\begin{theorem} \label{thm: gradient of dual objective}
    The function~$\omega \colon \setR^m \to \setR$ defined by~\cref{eq: omega} is continuously differentiable at every~$\lambda \in \setR^m$ and
    \begin{multline}
            \nabla \omega(\lambda) = g\left( \prox_{\frac{1}{\ell} \sum\limits_{i = 1}^{m} \lambda_i g_i} \left( y - \frac{1}{\ell} \sum_{i = 1}^{m} \lambda_i \nabla f_i(y) \right) \right) \\
                                                  + J_f(y) \left( \prox_{\frac{1}{\ell} \sum\limits_{i = 1}^{m} \lambda_i g_i} \left( y - \frac{1}{\ell} \sum_{i = 1}^{m} \lambda_i \nabla f_i(y) \right) - y \right) + f(y) - F(x)
    ,\end{multline}
    where~$\prox$ is the proximal operator~\cref{eq: proximal operator}, and~$J_f(y)$ is the Jacobian matrix at~$y$ given by
    \[
        J_f(y) \coloneqq \left( \nabla f_1(y), \dots, \nabla f_m(y) \right)^\top
    .\] 
\end{theorem}
\begin{proof}
    Define
    \[
        h(z, \lambda) \coloneqq \sum_{i = 1}^{m} \lambda_i g_i(z) + \frac{\ell}{2} \norm*{z - y + \frac{1}{\ell} \sum_{i = 1}^{m} \lambda_i \nabla f_i(y)}^2
    .\] 
    Clearly,~$h$ is continuous on~$\setR^n \times \setR^m$.
    Moreover,~$h_z(\cdot) \coloneqq h(z, \cdot)$ is continuously differentiable and
    \[
        \nabla_\lambda h_z(\lambda) = g(z) + J_f(y) \left( z - y + \frac{1}{\ell} \sum_{i = 1}^{m} \lambda_i \nabla f_i(y) \right) 
    .\] 
    Furthermore,
    \[
        \prox_{\frac{1}{\ell} \sum\limits_{i = 1}^{m} \lambda_i g_i}\left( y - \frac{1}{\ell} \sum_{i = 1}^{m} \lambda_i \nabla f_i(y) \right) = \argmin_{z \in \setR^n} h(z, \lambda)
    \] 
    is also continuous at every~$\lambda \in \setR^m$ (cf.~\cite[Theorem~2.26 and Exercise~7.38]{Rockafellar1998}).
    Therefore, the well-known result in first order differentiability analysis of the optimal value function~\cite[Theorem~4.13]{Bonnans2000} gives
    \begin{align}
            \MoveEqLeft \nabla_\lambda \left[ \ell \envelope_{\frac{1}{\ell} \sum\limits_{i = 1}^{m} \lambda_i g_i}\left( y - \frac{1}{\ell} \sum_{i = 1}^{m} \lambda_i \nabla f_i(y) \right) \right] \\
            ={}& g\left( \prox_{\frac{1}{\ell} \sum\limits_{i = 1}^{m} \lambda_i g_i} \left( y - \frac{1}{\ell} \sum_{i = 1}^{m} \lambda_i \nabla f_i(y) \right) \right) \\
               &+ J_f(y) \left( \prox_{\frac{1}{\ell} \sum\limits_{i = 1}^{m} \lambda_i g_i} \left( y - \frac{1}{\ell} \sum_{i = 1}^{m} \lambda_i \nabla f_i(y) \right) - y + \frac{1}{\ell} \sum_{i = 1}^{m} \lambda_i \nabla f_i(y) \right) 
    .\end{align}
    On the other hand, we have
    \begin{align}
            \MoveEqLeft \nabla_\lambda \left[ - \frac{1}{2 \ell} \norm*{\sum_{i = 1}^{m} \lambda_i \nabla f_i(y)}^2 + \sum_{i = 1}^{m} \lambda_i \left\{ f_i(y) - F_i(x) \right\} \right] \\
        &= - \frac{1}{\ell} J_f(y) \sum_{i = 1}^{m} \lambda_i \nabla f_i(y) + f(y) - F(x)
    .\end{align}
    Adding the above two equalities, we get the desired result.
\end{proof}
This theorem shows that the dual problem~\cref{eq: dual} is an $m$-dimensional differentiable convex optimization problem.
Hence, if we can compute the proximal operator of~$\sum_{i = 1}^{m} \lambda_i g_i$ quickly, then we can solve~\cref{eq: dual} using convex optimization techniques such as the interior point method~\cite{Boyd2004}.
In addition, for cases where $n \gg m$, the computational cost is much lower than solving the subproblem~\cref{eq: acc prox subprob} directly.
In particular, when~$m = 2$, eliminating a variable with~$\lambda_2 = 1 - \lambda_1$ reduces~\cref{eq: dual} to a one-dimensional optimization that can be solved quickly using, for example, Brent's method~\cite{Brent1973}.
Note, for example, that if~$g_i(x)=g_1(x)$ for all $i = 1, \dots, m$, or if~$g_i(x)=g_i(x_{I_i})$ and the index sets~$I_i$ do not overlap each other, then we can evaluate the proximal operator of~$\sum_{i = 1}^{m} \lambda_i g_i$ from the proximal operator of each~$g_i$.
Furthermore, even if there is an overlap, we can compute such a proximal operator immediately for special functions, for example, $m=2, g_1(x)=\norm*{x}_1,g_2(x)=\norm*{x}_2^2$ ($\lambda_1 g_1(x) + \lambda_2 g_2(x)$ is the elastic net~\cite{Zou2005} when~$\lambda_1 > 0$ and~$\lambda_2 > 0$. The elastic net has a proximal operator in closed-form~\cite[Section 6.5.3]{Parikh2014}).

\ifmain
\else
\fi

\section{Numerical experiments} \label{sec: experiments}
This section illustrates the proposed method's performance compared to the proximal gradient method without acceleration (\cref{alg: pgm}), and the algorithm below.
\begin{algorithm}[hbtp]
    \caption{Accelerated proximal gradient method for multiobjective optimization (without $f_i(y) - F_i(x)$)}
    \label{alg: acc-pgm-deprecated}
    \begin{algorithmic}[1]
        \Require Set~$x^0 = y^1 \in \dom F, \ell \ge L, \varepsilon > 0$.
        \Ensure $x^\ast$: A weakly Pareto optimal point
        \State $k \gets 1$
        \State $t_1 \gets 1$ \label{line: t ini deprecated}
        \Loop
        \State $x^k \gets$ the optimal solution of~\cref{eq: analogy} \label{line: x deprecated}
        \If{$\norm*{x^k - y^k}_\infty < \varepsilon$}
        \State \Return $x^k$
        \EndIf
        \State $t_{k + 1} \gets \sqrt{t_k^2 + 1 / 4} + 1/2$ \label{line: t rr deprecated}
        \State $\gamma_k \gets (t_k - 1) / t_{k + 1}$ \label{line: gamma deprecated}
        \State $y^{k + 1} \gets x^k + \gamma_k (x^k - x^{k - 1})$ \label{line: y deprecated}
        \State $k \gets k + 1$
        \EndLoop
    \end{algorithmic}
\end{algorithm}
Unlike the proposed \cref{alg: acc-pgm}, \cref{alg: acc-pgm-deprecated} does not include the term~$f_i(y) - F_i(x)$, which was the key to the proof of \cref{thm: acc conv rate}, in the subproblem solved in Step~\ref{line: x deprecated}.
Therefore, the convergence rate of \cref{alg: acc-pgm-deprecated} is still theoretically unknown.
However, since it is the easiest algorithm to conceive from the scalar optimization FISTA, and \cref{alg: acc-pgm-deprecated} is consistent with~\cite{ElMoudden2021} when~$g = 0$, we use it as a comparison in the numerical experiments.

\subsection{Test problems}
We generate a new list of convex multiobjective optimization test problems by processing the problem list of~\cite{Mita2019} based on the following three criteria:
\begin{itemize}
    \item Extracting convex problems: Since our proposed method is designed for convex problems, we selected only the convex problems from the original problem list.
    \item Dealing with various dimensions of~$n$: For some problems, we enhance the variety by using different values of~$n$.
    \item Including $g$: The original test problems include both constrained and unconstrained problems. For constrained problems, we set~$g_i$ as the indicator function~\cref{eq: indicator} corresponding to the constraint for every~$i = 1, \dots, m$. For unconstrained problems, we consider two types: $g_i = 0$ and $g_i = \norm{x - i + 1}_1 / [(i - 1) n]$ (i.e. using $\ell_1$-norm) for each~$i = 1, \dots, m$.
\end{itemize}
Based on this, the new list is given in~\cref{tab: test_problems}. 
\begin{table}[htbp]
  \centering
  \caption{List of test problems}
  \label{tab: test_problems}
  \begin{tabular}{lccc}
    \toprule
    Problem name & $m$ & $n$ & $g_i$ \\
    \midrule
    JOS1~\cite{Jin2001}     & $2$ & $5, 10, 20, 50, 100, 200, 500, 1000$ & $0$ or $\ell_1$ \\
    ZDT1~\cite{Zitzler2000} & $2$ & $5, 10, 20, 50, 100, 200, 500, 1000$ & indicator function \\
    SD~\cite{Stadler1992}   & $2$ & $4$ & indicator function \\
    TOI4~\cite{Toint1983}   & $2$ & $4$ & $0$ or $\ell_1$ \\
    TRIDIA~\cite{Toint1983} & $3$ & $3$ & $0$ or $\ell_1$ \\
    FDS~\cite{Fliege2009}   & $3$ & $5, 10, 20, 50, 100$ & $0$ or $\ell_1$ \\
    LFR1~\cite{More1981}    & $4$ & $30, 100, 1000$ & $0$ or $\ell_1$ \\
    \bottomrule
  \end{tabular}
\end{table}

\subsection{Experimental settings}
\label{sec: settings}
The experiments are carried out on a machine with $2.4$~GHz Intel Xeon Silver 4210R CPU and $64$~GB memory, implementing all codes in Python~3.9.5.
In all algorithms, we convert the subproblem into its dual as discussed in \cref{sec: subproblem} and solve it using the trust-region interior point method~\cite{byrd1999} with the scientific library SciPy. The stopping tolerance for solving the subproblem is $10^{-11}$, except for the difficult problem TRIDIA where we use $10^{-6}$. 
Also, we use backtracking procedure to determine a parameter~$\ell$, where the initial value of~$\ell$ is~$1$ and the constant multiplied to~$\ell$ is~$2$.
We set the general stopping criteria as~$\varepsilon = {10}^{-5}$ for each experiment.
Moreover, we choose~100 initial points, commonly for both algorithms, uniformly, and randomly between the bounds given in~\cite{Mita2019}.
The source code used here is available at \url{https://github.com/zalgo3/zfista}.

\subsection{Evaluation metrics}
We use the following metrics to assess the algorithms' performance:
\begin{itemize}
    \item \textbf{The number of iterations}: The number of iterations required to satisfy the stopping criteria.
    \item \textbf{Time}: The time needed to meet the stopping criteria.
    \item \textbf{Purity~\cite{Bandyopadhyay2004}}: The ratio of the solutions obtained by a given solver within the approximated Pareto frontier. Let~$PF_{p, s}$ be the set of function values of the solutions obtained by solver $s \in \mathcal{S}$ for problem $p \in \mathcal{P}$ that are not dominated by other solutions, and let~$PF_p$ be the set of~$\bigcup_{s \in \mathcal{S}} PF_{p, s}$ that are not dominated by other solutions. The purity is defined by~$\abs{PF_{p, s} \cap PF_p} / \abs{PF_p}$.
    \item \textbf{Hypervolume~\cite{Zitzler1999}}: The sum of the volumes of the hyperrectangles where the line segment connecting the reference point and each point of~$PF_{p, s}$ forms a diagonal. We set as the reference point the maximum value of each objective function in~$PF_p$.
    \item \textbf{Spread metrics ($\Gamma$ and~$\Delta$)~\cite{Custodio2011}}: The metric representing how well-distributed the obtained Pareto frontier is. Let $PF_{p, s} \cap PF_p$ be formed by~$F^1, \dots, F^N$. Assume that~$F_j^{a_j^i} \le F_j^{a_j^{i + 1}}$ for some~$\{a^i_j\} \subseteq \{1, \dots, m\}$ and for each~$i = 1, \dots, N, j = 1, \dots, m$. Moreover, set~$F_j^{a_j^0}$ and~$F_j^{a_j^{N + 1}}$ as the points in~$PF_p$ where~$F_j$ is largest and smallest, respectively. When~$N \ge 2$, the spread metrics~$\Gamma$ and~$\Delta$ are defined by
        \[
            \Gamma_{p, s} \coloneqq \max_{j = 1, \dots, m} \max_{i = 0, \dots, N} \delta_{i, j}
        \]
        and
        \[
            \Delta_{p, s} \coloneqq \max_{j = 1, \dots, m} \frac{\delta_{0, j} + \delta_{N, j} + \sum_{i = 1}^N \abs{\delta_{i, j} - \bar{\delta}_j}}{\delta_{0, j} + \delta_{N, j} + (N - 1) \bar{\delta}_j},
        \]
        where~$\delta_{i, j} \coloneqq F_j^{a_j^{i + 1}} - F_j^{a_j^i}$ and~$\bar{\delta}_j \coloneqq \sum_{i = 1}^N \delta_{i, j} / N$.
        On the other hand when~$N \le 1$, we define~$\Gamma_{p, s} = \Delta_{p, s} = \infty$.
\end{itemize}

We also obtained performance profiles~\cite{Dolan2002} for each of the evaluation metrics to provide a comprehensive comparison of the algorithms.
Suppose that a metric~$t_{p, s}$ is defined for a solver~$s \in \mathcal{S}$ and a problem~$p \in \mathcal{P}$.
We assume that the smaller~$t_{p, s}$ is, the better.
The performance profile~$R_s(\tau)$ of a solver~$s \in \mathcal{S}$ is defined as
\[
    R_s(\tau) \coloneqq \frac{1}{\abs{\mathcal{P}}} \abs*{p \in \mathcal{P} \mid r_{p, s} \le \tau},
\]
where~$r_{p, s}$ is the performance ratio given by~$r_{p, s} \coloneqq t_{p, s} / \min_{s \in \mathcal{S}} t_{p, s}$.
Note that for hypervolume and spread metrics, we took the reciprocal when calculating the performance ratio, as larger metric values correspond to better performance for them.

\subsection{Results of the experiments}

\begin{figure}[htbp]
    \centering
    \begin{minipage}[b]{.48\linewidth}
        \centering
        \includegraphics[clip, width=\linewidth]{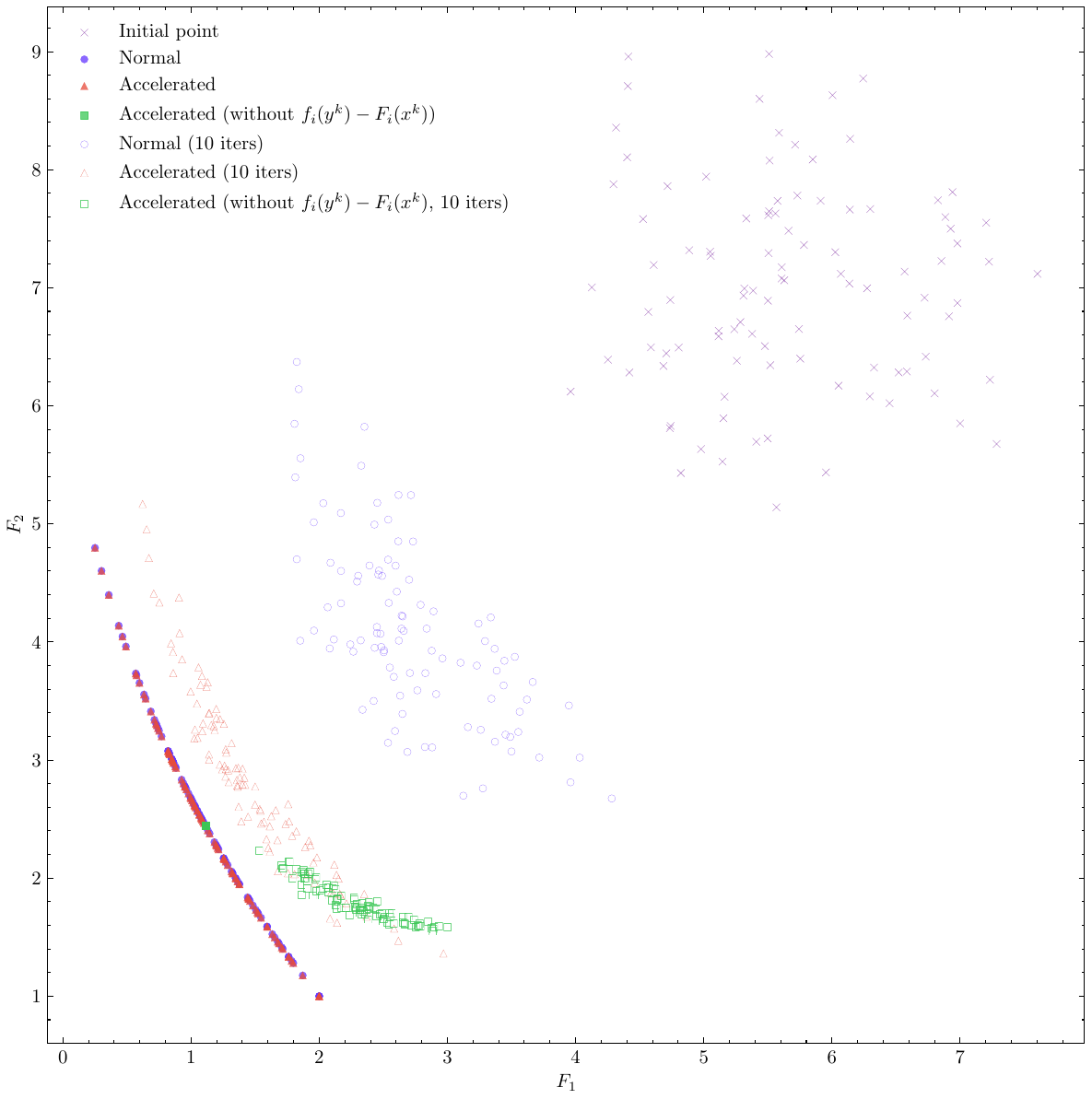}
        \caption{Objective function values for problem~JOS1 with $n=50$, and $\ell_1$~norm for $g_i$}
        \label{fig: JOS1_n_50_l1_objective}
    \end{minipage}
    \hfill
    \begin{minipage}[b]{.5\linewidth}
        \centering
        \includegraphics[clip, width=\linewidth]{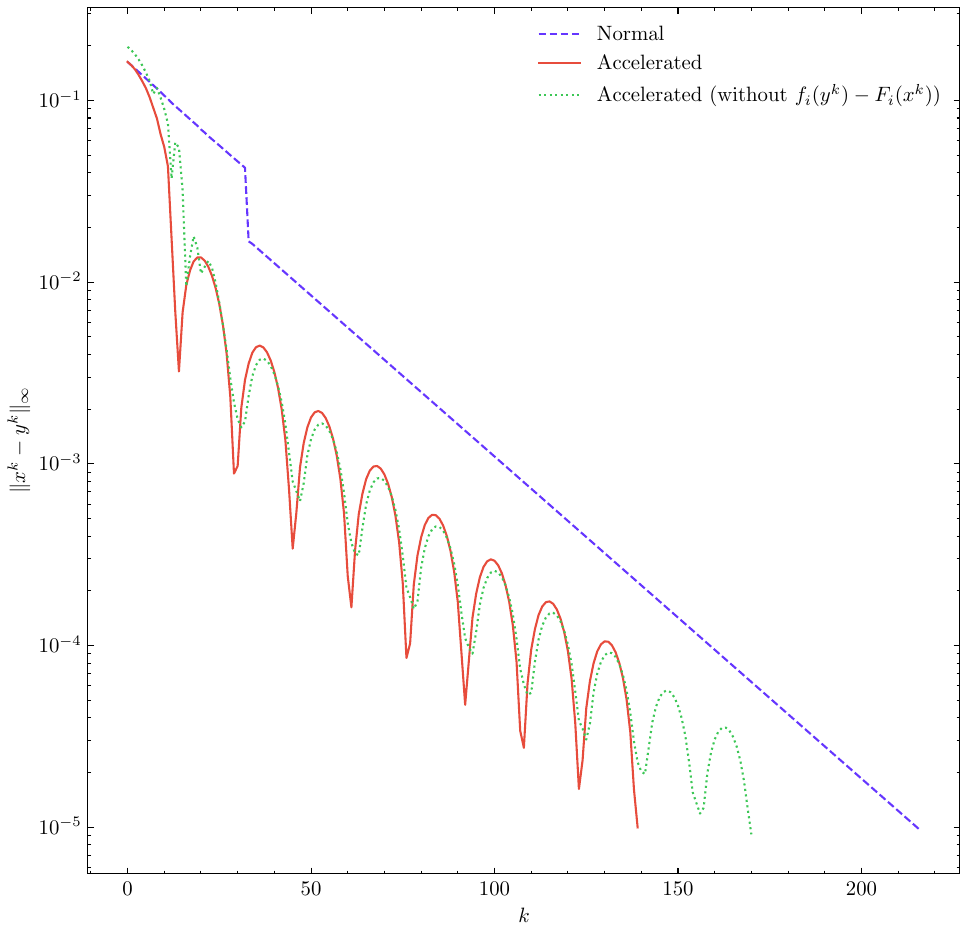}
        \caption{An example of~$\setof*{\norm*{x^k - y^k}_\infty}$ for problem~JOS1 with $n=50$, and $\ell_1$~norm for $g_i$}
        \label{fig: JOS1_n_50_l1_error}
    \end{minipage}
\end{figure}

Let us first illustrate the behaviour of the algorithms. For this, we take the problem JOS1~\cite{Jin2001} with $n = 50$ and $g_i$ as the $\ell_1$-norm. In \Cref{fig: JOS1_n_50_l1_objective}, we plot the objective function values for $k = 0$ (i.e., at the initial points), $k = 10$, and the terminal points of each algorithm, respectively. The set of terminal points are in fact the Pareto solutions obtained. Here, ``Normal'', ``Accelerated'', and ``Accelerated (without $f_i(y) - F_i(x)$)'' means, respectively, \cref{alg: pgm}, \cref{alg: acc-pgm} and \cref{alg: acc-pgm-deprecated}.
As we can see, all the algorithms were able to find a wide range of Pareto solutions in this case. However, the objective function values at~$k = 10$ are smaller when using the accelerated \cref{alg: acc-pgm} and \cref{alg: acc-pgm-deprecated}. Moreover, from \cref{fig: JOS1_n_50_l1_error}, we see that \cref{alg: acc-pgm} and \cref{alg: acc-pgm-deprecated} converge faster than the non-accelerated \cref{alg: pgm}, despite oscillations. In this example, we can also see that \cref{alg: acc-pgm} were faster and obtained a more uniform Pareto frontier than \cref{alg: acc-pgm-deprecated}.

\begin{figure}[htbp]
    \centering
    \begin{minipage}[b]{.49\linewidth}
        \centering
        \includegraphics[clip, width=\linewidth]{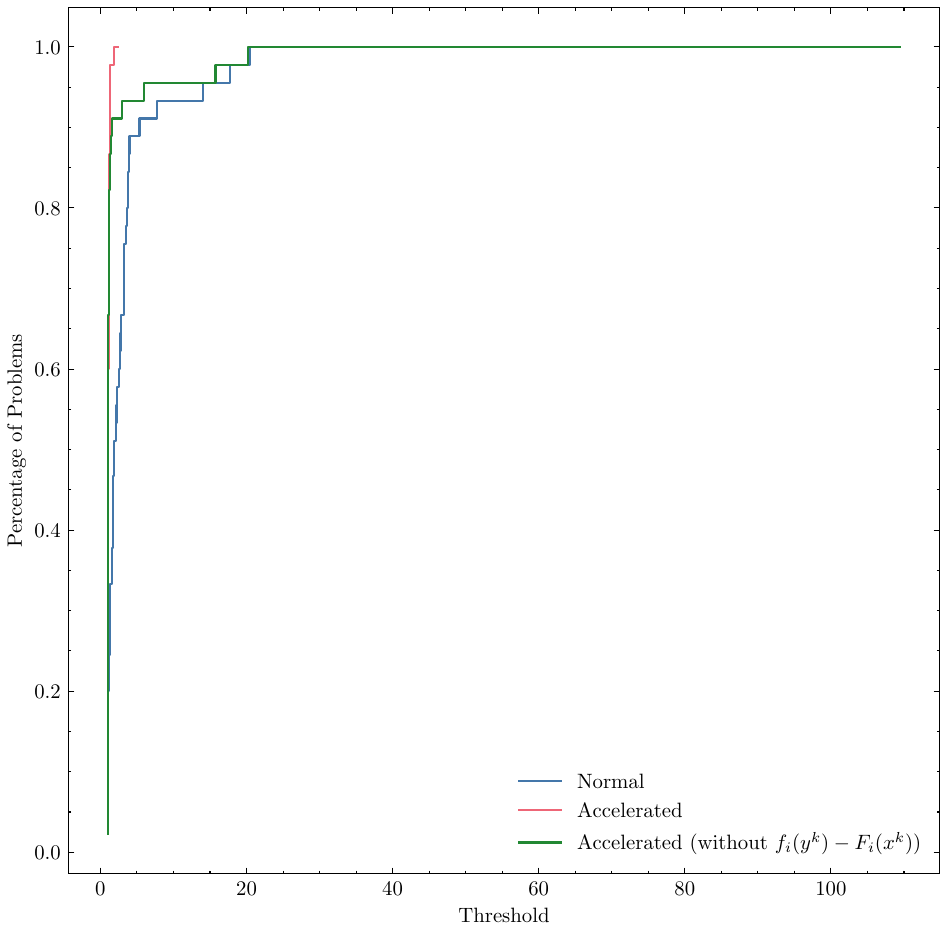}
        \caption{Performance profile:\\ computational time}
        \label{fig: pf_computational_time}
    \end{minipage}
    \hfill
    \begin{minipage}[b]{.49\linewidth}
        \centering
        \includegraphics[clip, width=\linewidth]{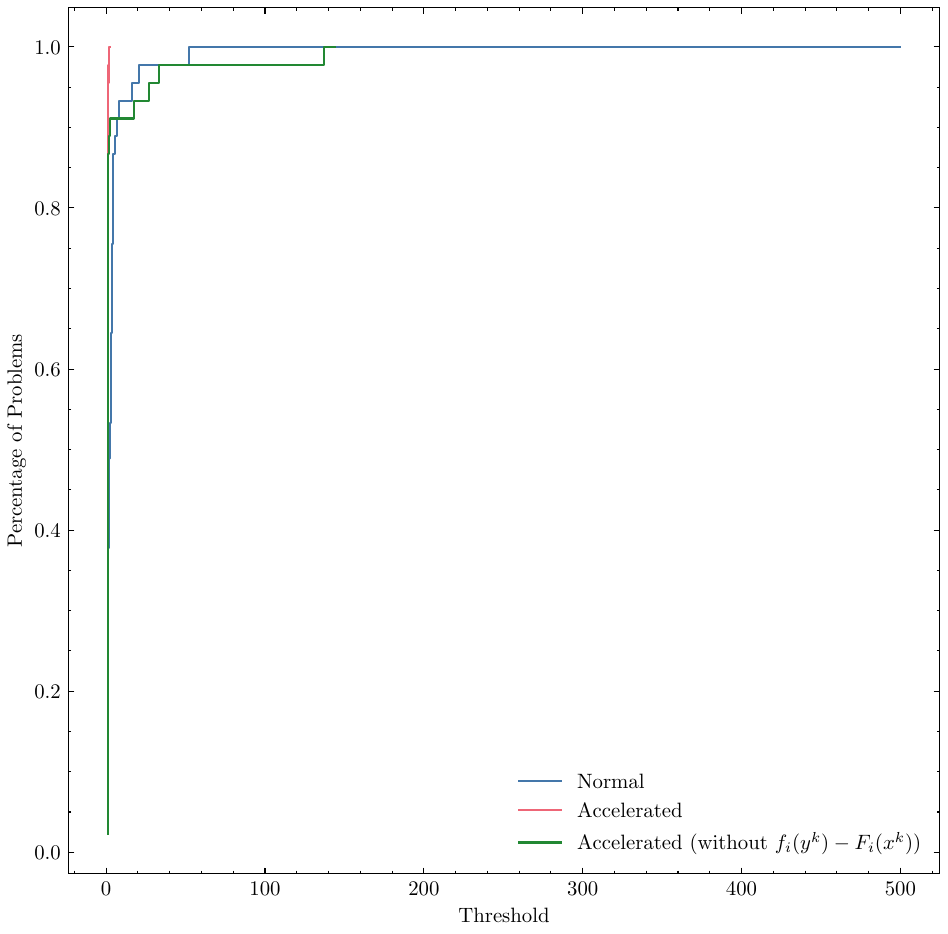}
        \caption{Performance profile:\\ iteration counts}
        \label{fig: pf_iteration_counts}
    \end{minipage}
\end{figure}

We now check the performance of the algorithms. As we explained in~\cref{sec: settings}, for each problem of~\cref{tab: test_problems}, we run the algorithms with $100$~different initial points. \cref{tab: computational_time} shows the average of the computational time and iteration counts for each problem. For problems with different values of $n$, we just show the smallest and the biggest $n$ for convenience. From the table, it is possible to see that acceleration is in general more efficient in terms of time. In fact, by checking the performance profiles given in~\cref{fig: pf_computational_time} and~\cref{fig: pf_iteration_counts}, we observe that our proposed~\cref{alg: acc-pgm} performs better in terms of time and iteration counts. It is interesting to see from \cref{tab: computational_time}, however, that there are cases where~\cref{alg: acc-pgm} does not perform well. 

\begin{table}[htbp]
  \centering
  \caption{Average computational costs}
  \label{tab: computational_time}
  \begin{tabular}{l@{\,}cc|rrr|rrr}
    \toprule
    \multirow{2}{*}{Problem} & \multirow{2}{*}{$n$} & \multirow{2}{*}{$g_i$}
    & \multicolumn{3}{c}{Total time (s)} & \multicolumn{3}{|c}{Iteration counts} \\
    & & & Alg.~\ref{alg: pgm} & Alg.~\ref{alg: acc-pgm} & Alg.~\ref{alg: acc-pgm-deprecated} &
    Alg.~\ref{alg: pgm} & Alg.~\ref{alg: acc-pgm} & Alg.~\ref{alg: acc-pgm-deprecated} \\
    \midrule
    JOS1 & $5$    &      $0$ & $0.032$    & $0.029$    & $0.027$    & $23.82$   & $27.89$  & $27.89$ \\
    JOS1 & $5$    & $\ell_1$ & $0.783$    & $0.355$    & $0.336$    & $22.20$   & $21.26$  & $28.09$ \\
    JOS1 & $1000$ &      $0$ & $3.674$    & $0.260$    & $0.207$    & $3203.05$ & $155.00$ & $155.00$ \\ 
    JOS1 & $1000$ & $\ell_1$ & $183.957$  & $47.913$   & $46.197$   & $2901.50$ & $732.72$ & $644.11$ \\
    ZDT1 & $5$    & ind.     & $0.743$    & $0.279$    & $0.234$    & $38.81$   & $11.03$  & $8.87$ \\
    ZDT1 & $1000$ & ind.     & $1.738$    & $0.840$    & $0.634$    & $32.68$   & $14.30$  & $9.986$ \\
    SD   & $4$    & ind.     & $1.063$    & $0.806$    & $1.026$    & $36.58$   & $33.02$  & $32.94$ \\
    TOI  & $4$    &      $0$ & $0.013$    & $0.018$    & $0.015$    & $3.95$    & $4.57$   & $5.18$ \\
    TOI  & $4$    & $\ell_1$ & $1.109$    & $0.841$    & $1.035$    & $20.95$   & $18.41$  & $22.90$ \\
    TRIDIA & $3$  &      $0$ & $94.447$   & $0.981$    & $5.842$    & $3177.21$ & $6.35$   & $172.89$ \\
    TRIDIA & $3$  & $\ell_1$ & $79.562$   & $3.892$    & $78.616$   & $1348.16$ & $25.80$  & $860.90$ \\    
    FDS  & $5$    &      $0$ & $22.897$   & $12.719$   & $14.433$   & $286.4$   & $152.35$ & $170.83$ \\
    FDS  & $5$    & $\ell_1$ & $16.611$   & $12.150$   & $1330.424$ & $127.48$  & $91.39$  & $13178.25$ \\
    FDS  & $100$  &      $0$ & $3805.058$ & $4007.926$ & $3607.842$ & $644.45$  & $117.27$ & $158.72$ \\
    FDS  & $100$  & $\ell_1$ & $4773.474$ & $5412.802$ & $5880.479$ & $767.81$  & $177.37$ & $474.55$ \\
    LFR1 & $30$   &      $0$ & $4.904$    & $8.335$    & $4.362$    & $9.18$    & $11.67$  & $6.69$ \\
    LFR1 & $30$   & $\ell_1$ & $10.337$   & $10.184$   & $160.399$  & $8.91$    & $11.4$   & $1224.37$ \\
    LFR1 & $1000$ &      $0$ & $10.928$   & $13.804$   & $10.916$   & $8.54$    & $10.07$  & $8.36$ \\
    LFR1 & $1000$ & $\ell_1$ & $26.669$   & $31.038$   & $26.566$  & $8.55$    & $10.31$   & $8.56$ \\ 
    \bottomrule 
  \end{tabular}
\end{table}

\begin{figure}[htbp]
    \centering
    \begin{minipage}[b]{.49\linewidth}
        \centering
        \includegraphics[clip, width=\linewidth]{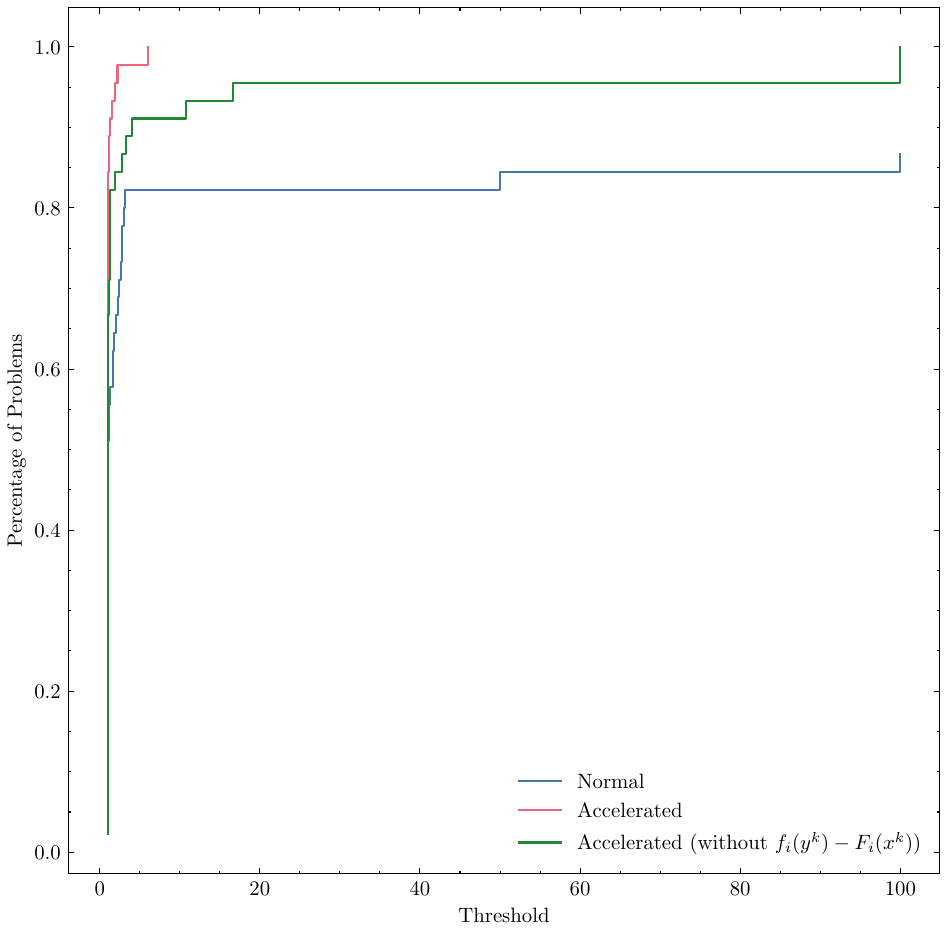}
        \caption{Performance profile:\\ purity}
        \label{fig: pf_purity}
    \end{minipage}
    \hfill
    \begin{minipage}[b]{.49\linewidth}
        \centering
        \includegraphics[clip, width=\linewidth]{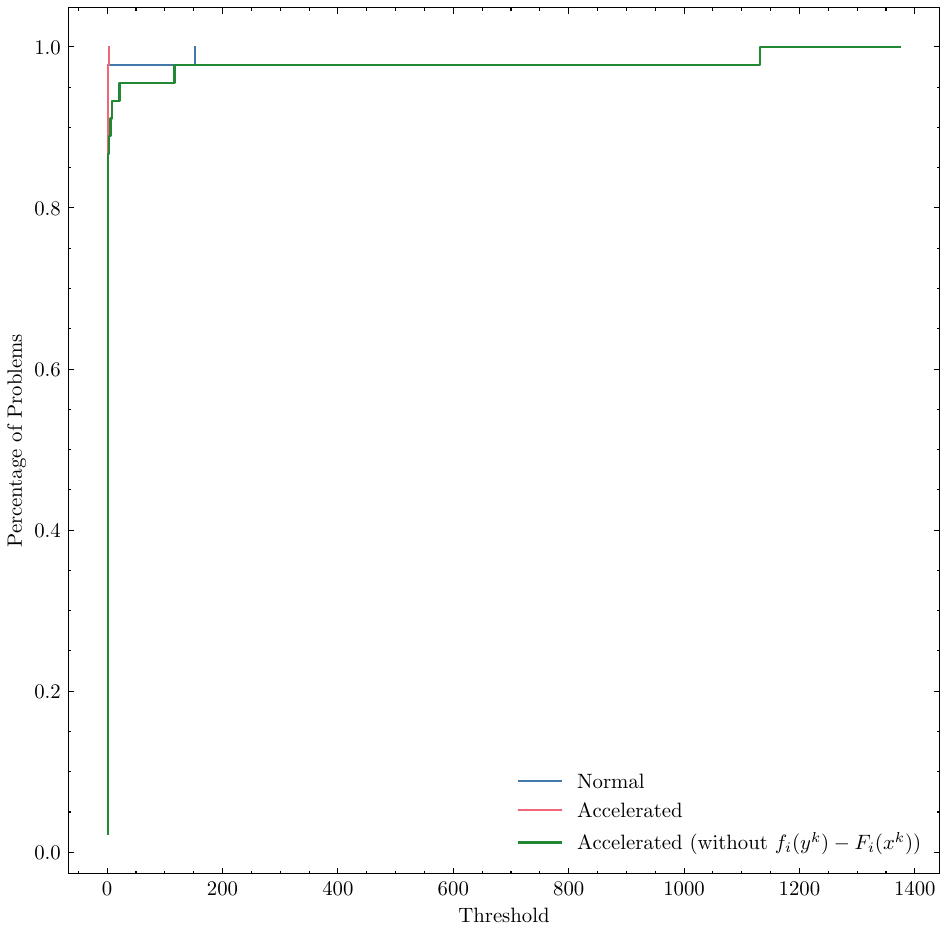}
        \caption{Performance profile: hypervolume}
        \label{fig: pf_hypervolume}
    \end{minipage}
\end{figure}

\begin{figure}[htbp]
    \centering
    \begin{minipage}[b]{.49\linewidth}
        \centering
        \includegraphics[clip, width=\linewidth]{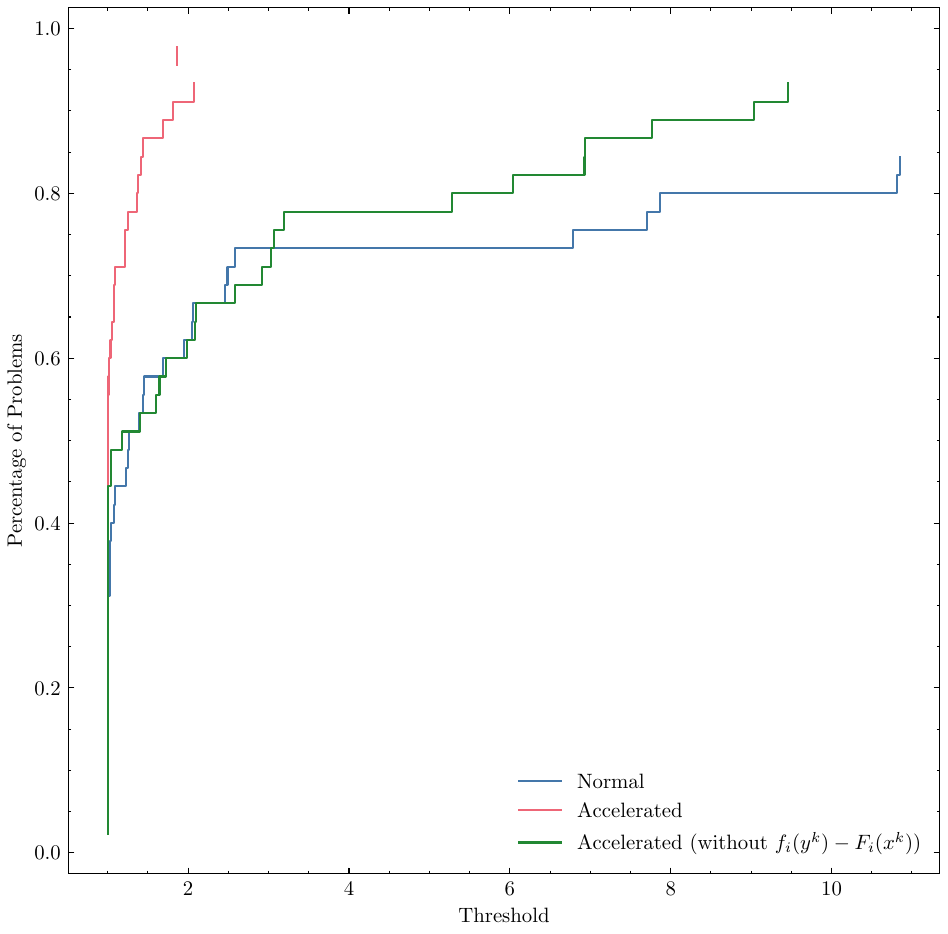}
        \caption{Performance profile:\\ spread metrics ($\Gamma$)}
        \label{fig: pf_gamma}
    \end{minipage}
    \hfill
    \begin{minipage}[b]{.49\linewidth}
        \centering
        \includegraphics[clip, width=\linewidth]{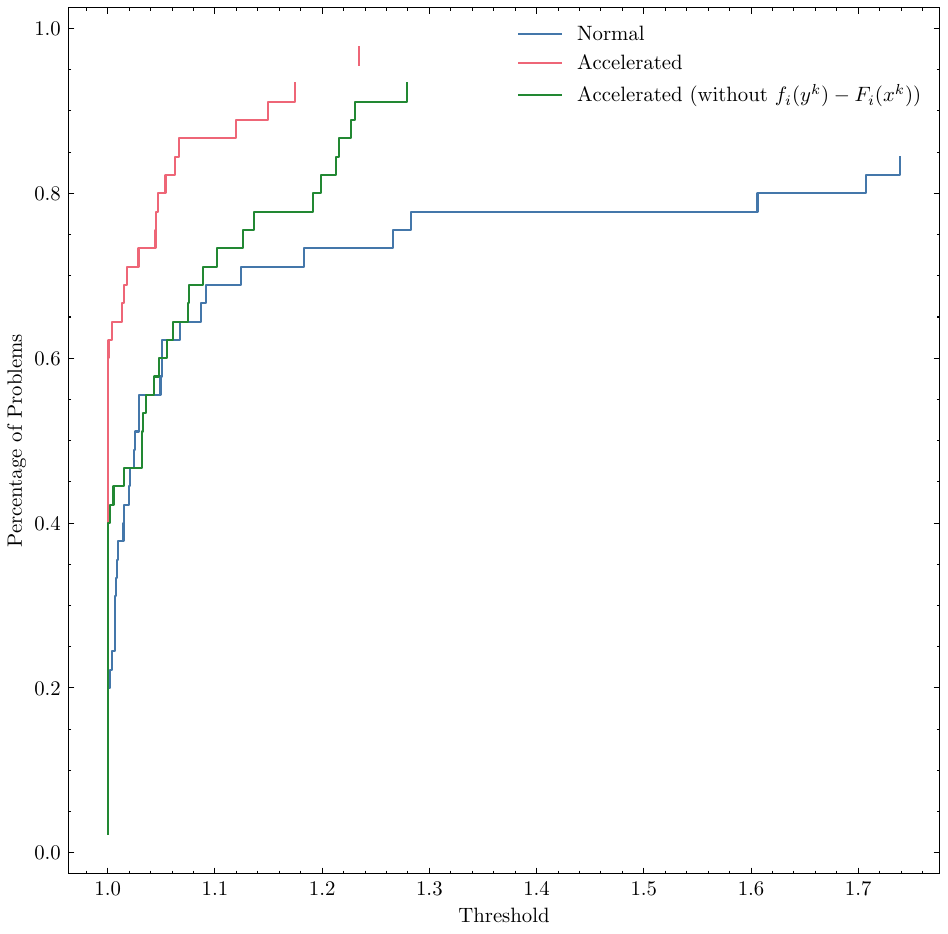}
        \caption{Performance profile:\\ spread metrics ($\Delta$)}
        \label{fig: pf_delta}
    \end{minipage}
\end{figure}

Besides the performance, it is usually important to see how good the Pareto frontier is. Thus, once again we show performance profiles, this time for purity (\cref{fig: pf_purity}), hypervolume (\cref{fig: pf_hypervolume}), spread metric $\Gamma$ (\cref{fig: pf_gamma}) and spread metric $\Delta$ (\cref{fig: pf_delta}). Clearly, our proposed~\cref{alg: acc-pgm} outperforms the other two algorithms, obtaining better Pareto frontiers. We can thus conclude that at least among the test problems considered, \cref{alg: acc-pgm} seem promising both in terms of performance and uniform Pareto frontiers.

\ifmain
\else
\fi

\section{Conclusion} \label{sec: conclusion}
By putting information of the previous points into the subproblem, we have successfully accelerated the proximal gradient method for multiobjective optimization and proved its convergence rate under natural assumptions, which was an open problem.
Moreover, we showed an efficient way of computing the subproblem via its dual.
As the experiments suggested, the proposed methods are also effective from the numerical point of view.

This paper shows the convergence rate for the sequence of the merit function values and the classical global convergence concerning accumulation points but does not provide the global convergence of the sequence of iterates itself.
For single-objective optimization, by changing the update rule for the parameter~$t_k$, Chambolle and Dossal have proposed a variant with the iterates' global convergence~\cite{Chambolle2015}.
In the multi-objective optimization problem, it may also be possible to modify the algorithm similarly and obtain the global convergence of iterates.
Moreover, since many schemes for single-objective optimization had been developed, following the idea of Nesterov's acceleration technique, this paper may also contribute to the development of various multiobjective optimization methods.
Extensions to vector optimization and its generalization, the vector optimization problem with variable domination structure~\cite{Bouza2022,Kobis2022}, may also be worth considering.
Such extensions will be subjects of future works.

\ifmain
\else
\fi

\ifpreprint
    \section*{Acknowledgments}
\fi
This work was supported by the Grant-in-Aid for Scientific Research (C) (21K11769 and 19K11840) and Grant-in-Aid for JSPS Fellows (20J21961) from the Japan Society for the Promotion of Science.

\bibliographystyle{jorsj}
\bibliography{library}

\end{document}

\bibliographystyle{junsrt}
\bibliography{}
\end{document}